\documentclass{amsart}
\usepackage{amssymb}
\newtheorem{lemma}{Lemma}[section]
\newtheorem{prop}[lemma]{Proposition}

\newtheorem{theo}[lemma]{Theorem}

\newtheorem{corol}[lemma]{Corollary}

\newcommand{\C}{{\mathbb C}}

\newcommand{\Z}{{\mathbb Z}}

\newcommand{\CC}{\hbox{{$\mathcal C$}}}
\newcommand{\CM}{\hbox{{$\mathcal M$}}}
\newcommand{\CN}{\hbox{{$\mathcal N$}}}
\newcommand{\CF}{\hbox{{$\mathcal F$}}}

\newcommand{\extd}{{d}}
  \newcommand{\del}{\partial}
\newcommand{\image}{{\rm image}}

\newcommand{\isom}{{\cong}}
\newcommand{\eps}{{\varepsilon}}
\newcommand{\tens}{\mathop{\otimes}}
\newcommand{\la}{{\triangleright}}
\newcommand{\ra}{{\triangleleft}}

\newcommand{\Ad}{{\rm Ad}}
\newcommand{\id}{{\rm id}}
\newcommand{\<}{\langle}
\renewcommand{\>}{\rangle}

\newcommand{\eproof}{$\quad \diamond$\bigskip}
\newcommand{\und}[1]{{\underline {#1}}}
\newcommand{\eqn}[2]{\begin{equation}#2\label{#1}\end{equation}}
\newcommand{\lcross}{{>\kern -4pt\triangleleft}}
\newcommand{\rcross}{{\triangleright \kern -4pt  <}}
\newcommand{\bla}{{\blacktriangleright}}
\newcommand{\bra}{{\blacktriangleleft}}
\newcommand{\rcrossco}{{\bla\kern -4pt <}}
\newcommand{\lrbicross}{{\bla \kern-1pt\triangleleft}}
\newcommand{\rlbicross}{{\triangleright\kern-1pt\bra}}

\newcommand{\tra}{{\tilde{\triangleleft}}}

\begin{document}

\title[CARTAN CALCULUS ON BICROSSPRODUCTS]{\rm \large
CARTAN CALCULUS FOR QUANTUM DIFFERENTIALS ON BICROSSPRODUCTS}
\author{F. Ngakeu, S. Majid, J-P. Ezin}%

\address{F.N. and J-P. E: Institut de Math\'ematiques et de
Sciences Physiques\\
BP 613 Porto-Novo, Benin}
\address{S.M: School of Mathematical Sciences\\
Queen Mary, University of London\\ 327 Mile End Rd,  London E1
4NS, UK.}

\thanks{F.N. acknowledges financial support from CIUF-CUD}
\thanks{S.M. is a Royal Society University Research Fellow at QMUL}

\date{5/2002; revised 6/2003}%

\maketitle

\begin{abstract}
We provide the Cartan calculus for bicovariant differential forms
on bicrossproduct quantum groups $k(M)\lrbicross kG$ associated to
finite group factorizations $X=GM$ and a field $k$. The
irreducible calculi are associated to certain conjugacy classes in
$X$ and representations of isotropy groups. We find the full
exterior algebras and show that they are inner by a bi-invariant
1-form $\theta$ which is a generator in the noncommutative de Rham
cohomology $H^1$. The special cases where one subgroup is normal
are analysed.  As an application, we study the noncommutative
cohomology on the quantum codouble
 $ D^*(S_3)\isom k(S_3)\lrbicross k\Z_6$ and the quantum double
 $D(S_3)=k(S_3)\lcross kS_3$, finding
respectively a natural calculus and a unique calculus with
$H^0=k.1$.
\end{abstract}

\section{Introduction}

There has been a lot of interest in recent years in finite groups
$M$, say, as noncommutative differential geometries (even though
the algebra of functions $k(M)$, $k$ a field, is commutative), see
\cite{ASitarz,KBresser,MAjRiema,MAjRai,MAjYangmills,ngaMAj}. The
bicovariant differential calculi on $k(M)$ are defined by
conjugacy classes $\CC\subset M$ not containing the group identity
and defined in practice by the Cartan calculus consisting of a
basis $\{e_a:\ a\in \CC\}$ of left-invariant differential 1-forms
and the bimodule and exterior derivative relations \eqn{dk(G)}{
\extd f=\sum_{a\in \CC} (R_a(f)-f)e_a,\quad e_a f =R_a(f)e_a,\quad
\forall f\in k(M)} where $R_a$ denotes right multiplication on the
group. It turns out in this way that there is an entire geometry
and Lie theory of finite groups. Another feature is that the
calculus is inner in the sense that there exists an element
$\theta=\sum_a e_a$ such that $\extd f=[\theta,f]$. Graded
commutator with $\theta$ similarly defines the differential in
higher degree, while a certain braiding $\Psi$ describes the
skew-symmetrization of basic 1-forms.

Since the suitable dual of a Hopf algebra is also a Hopf algebra,
one has another class of models where the `coordinate' algebra is
the group algebra $kG$, say, for a finite group $G$. If $G$ is
nonAbelian this is now genuinely noncommutative. Such objects
provide the first examples of noncommutative geometry which is
strictly noncommutative in both the quantum groups approach and
the Connes and operator theory approach (as for example in the
Baum-Connes theory for the K-theory of $\C G$ in terms of $EG$
\cite{BaumCon}). Differential calculi in this case were classified
in \cite{MAj1} and are given by irreducible right-representations
$V$ and vectors $\theta\in V/k$ (only the class of $\theta\in V$
controls the calculus). Here the invariant 1-forms are labelled by
a basis $e_i\in V$ and the calculus has the form \eqn{dkG}{ \extd
u=u \theta*(u-1),\quad e_i u=u (e_i*u),\quad \forall u\in G} where
$*$ denotes the right action. The calculus is inner via the chosen
$\theta$. Such models in the Lie setting would be the Hopf algebra
$U(g)$ where $g$ is a Lie algebra, for example $U(su_2)$ leads to
the `fuzzy sphere'. The Abelian discrete group case is also useful
e.g. after twisting to describe Clifford algebras as
noncommutative spaces and to describe noncommutative tori at the
algebraic level.

In the present paper we extend the above formulae to the next more
complicated finite noncommutative geometry in this family, namely
to bicrossproduct quantum groups \cite{Kac,Tak,MAjphy}
$k(M)\lrbicross kG$ where the above two models are `smashed
together'. These are now genuine noncommutative and
noncocommutative quantum groups. They also have the self-dual-type
feature namely the dual is $kM\rlbicross k(G)$ of the same
bicrossproduct type. For Lie groups they were proposed as
nontrivial noncommutative geometries (in connection with quantum
gravity) in \cite{MAjPla}  and as quantum Poincar\'e groups of
noncommutative spacetimes in \cite{MaRue}. More recently they have
played a role in computing cyclic cohomology\cite{ConMos} as well
as in the renormalisation of quantum field theories\cite{ConKre}.
Here they play a role linked to diffeomorphism invariance. The
finite group case is intimately linked to set-theoretic solutions
of the Yang-Baxter equations and over $\C$ was characterised by Lu
as Hopf algebras with positive basis\cite{Lu}. For all these
reasons it is clear that such bicrossproduct quantum groups should
be an important next most complicated and truly `quantum' source
of examples after the finite group cases. Their noncommutative
differential geometry, however, is very little explored and
explicit formulae for their differential structure, a prerequisite
for any actual computations and applications of the geometry, have
been totally lacking. We provide these now, in Section~3. Sections
4,5 cover special semidirect cases where either $M$ or $G$ are
normal.   The semidirect case in Section~4 also includes the
important case of the quantum codouble $D^*(G)=k(G)\rcrossco kG$
of a finite group, where we find a natural calculus induced from
one on $k(G)$ defined by a conjugacy class in $G$. Section~6
applies our Cartan calculus to explicit computations of
noncommutative de Rham cohomology, which turns out to be
nontrivial. The noncommutative differential geometry of the
quantum double $D(S_3)$ in Section~6.3, particularly, should be
physically interesting in connection with finite conformal field
theory and finite versions of fuzzy spheres. We find a unique
calculus with the connectedness property $H^0=k$.

Our starting point, in the preliminary Section~2, is the known but
nonconstructive classification theorem \cite{BM} for bicovariant
differentials on bicrossproducts due to E. Beggs and one of the
present authors. From the Woronowicz theorem\cite{W} one knows
that bicovariant calculi are classified by Ad-stable right ideals
in the augmentation ideal of the Hopf algebra. It was shown in
\cite{BM} that these are in 1-1 correspondence with certain
equivalence classes in the group $X=GM$ which determines the
bicrossproduct. We recall that if $X$ is a group factorization (in
the sense of two subgroups $G,M$ such that the product $G\times
M\to X$ is bijective) then each group acts on the other by actions
$\la,\ra$ defined by $su=(s\la u)(s\ra u)$ for $u\in G$ and $s\in
M$. They obey
\begin{eqnarray}\label{matchingcond}
s\ra e&=&s,\quad e\la u=u,\quad s\la e=e,\quad e\ra u=e \nonumber\\
(s\ra u)\ra v&=&s\ra (uv),\quad s\la(t\la u)=(st)\la u \\
s\la (uv)&=&(s\la u)((s\ra u)\la v),\quad (st)\ra u=(s\ra (t\la
u))(t\ra u) \nonumber
\end{eqnarray}
and conversely such a matched pair of actions allows to
reconstruct $X=G\bowtie M$ by a double cross product
construction\cite{MAjFounda}. Moreover, at least in the finite
case it means that the group algebra $kG$ acts on $k(M)$ and
$k(M)$ coacts on $kG$. The bicrossproduct $k(M)\lrbicross kG$ is
by definition the cross product algebra $\lcross$ by the action
and cross coproduct coalgebra $\rcrossco$ by the coaction.
Section~2 recalls the Beggs-Majid result with a slightly more
explicit description as the decomposition into conjugacy classes
of a certain $Z\subset X$. We also make a shift of conventions
from left modules to right modules which is not straightforward.
Our goal from this starting point is then to find a suitable basis
for the invariant differential forms and the Cartan calculus for
the differential structure. We find
(Theorem~\ref{cartanrelations}) that there is indeed a natural
choice of such a basis $\{e_a\}$ dual to a basis $\{f_a\}$ of the
quantum tangent space $L$, identified with a subrepresentation
under an action of $D(X)$ on $kX$. Hence there is the induced
$X$-graduation  $||.\;||$ on $L$ which factorizes as
$$||f_a||=\<f_a\>^{-1}|f_a|,$$
 say, in $MG$, and an induced right action $*$ of $X$ on $\{e_a\}$.
 Then (Theorem.~\ref{cartanrelations})
\[ e_{a}f=R_{\<f_a\>}(f)e_{a},\quad e_{a}u=(\<f_a\>\la
u)e_{a}*u\]
\[ \extd
f=\sum_{a}c_a(R_{\<f_a\>}(f)-f)e_{a},\quad \extd
u=\sum_{a}c_a((\<f_a\>\la u)e_{a}*u-ue_a)\]  where
$c_a=<\delta_{\<f_a\>},f_{a}>$ is defined by the pairing between
$kX$ and $k(X)$. We also find that the calculus is again inner.
These structures, and $\theta=\sum_{a}c_ae_a$ `unify' the two
extreme cases above when either $G$ or $M$ is trivial. Note that
there is no `algorithm' from \cite{BM} leading from the
classification to a suitable basis and resulting Cartan calculus
needed for practical applications, so that the work in the present
sequel is required. Further new results are the inner property and
that $\theta$ is a generator of the noncommutative de Rham
cohomology.

\subsection*{Preliminaries}

Here we collect all the basic definitions needed in the paper. We
work over a field $k$ of characteristic zero. Let $X=GM$ be a
finite group factorization. The bicrossproduct Hopf algebra
$A=k(M)\lrbicross kG$ has basis $\delta_s \otimes u$ where $s\in
M,u\in G$ and $\delta_s$ is the Kronecker delta-function in
$k(M)$. The product, coproduct $\Delta:A\to A\tens A$, counit
$\eps:A\to k$ and `coinverse' or antipode $S:A\to A$  for a Hopf
algebra are
 \eqn{relationsofAa}{ (\delta_s \otimes u)(\delta_t
\otimes v)=\delta_{s\ra u,t}(\delta_s \otimes uv),\quad  \Delta
(\delta_s \otimes u)=\sum\limits_{ab=s}\delta_a \otimes b\la
u\otimes \delta_b \otimes u} \eqn{relationsofAb}{ 1=\sum\limits_s
\delta_s \otimes e,\quad \varepsilon (\delta_s \otimes
u)=\delta_{s,e},\quad S(\delta_s \otimes u)=\delta_{(s\ra
u)^{-1}}\otimes (s\la u)^{-1}.} We use here the conventions and
notations for Hopf algebras in \cite{MAjFounda}. Thus,
$\Delta,\eps$ are algebra maps and coassociative (they define an
algebra on the dual) and $S$ obeys $\sum
(Sa_{(1)})a_{(2)}=\eps(a)1=\sum a_{(1)}(Sa_{(2)})$ for all $a$ if
we use the `Sweedler notation' $\Delta a=\sum a_{(1)}\tens
a_{(2)}$. The point of view in the paper is that $A$ is {\em like}
functions on a group and $\Delta,\eps,S$ encode the `group'
structure. Similarly, an action of this `group' is expressed as a
coaction of $A$, which is like an action but with arrows reversed.
Meanwhile, the dual $H=A^*=kM\rlbicross k(G)$ is also a
bicrossproduct, with \eqn{relationsofHa}{ (s\otimes
\delta_u)(t\otimes \delta_v)=\delta_{u,t\la v}(st\otimes
\delta_v),\quad \Delta(s\otimes \delta_u)=
\sum\limits_{xy=u}s\otimes \delta_x\otimes s\ra x\otimes \delta_y}
\eqn{relationsofHb}{ 1=\sum\limits_u e\otimes \delta_u,\quad
\varepsilon(s\otimes \delta_u) =\delta_{u,e},\quad  S(s\otimes
\delta_u)=(s\ra u)^{-1}\otimes \delta_{(s\la u)^{-1}}}

We use the Drinfeld quantum double $D(H)=H^{*op}\bowtie H$ built
on $H^{*}\otimes H$ in the double cross product form\cite{MAjphy},
see \cite{MAjFounda}. In the present case of $H=kM\rlbicross
k(G)$, the double was computed in \cite{BGM} and the cross
relations between $H$ and $H^{*op}$ are
\begin{eqnarray}
(1\otimes t\otimes \delta _v)(\delta_s \otimes u\otimes 1)=
\delta_{s'}\otimes u'\otimes t' \otimes \delta_{v'}
\end{eqnarray}
where
\begin{eqnarray}\label{theprimes}
s '&=&(t\ra (s\la u)^{-1})s(t\ra vu^{-1})^{-1}, \quad u'=(t\ra
vu^{-1})\la u \\ \nonumber t'&=&t\ra (s\la u)^{-1}, \quad v'=(s\la
u)vu^{-1}
\end{eqnarray}
obeying
\begin{eqnarray*}
t'\ra v'&=&t\ra vu^{-1}, \quad s' \la u'
=(t\la (s\la u)^{-1})^{-1} \nonumber\\
t'\la v'&=& (s' \la u')(t\la vu^{-1}) ,\quad s' \ra u'=t(s\ra
u)(t\ra v)^{-1}
\end{eqnarray*}
We use, and will freely use basic identities such as:
\begin{eqnarray*}
t^{-1}\ra (t\la u)&=&(t\ra u)^{-1},\quad (t\ra u)\la u^{-1}=(t\la
u)^{-1}
\end{eqnarray*}
\begin{eqnarray}\label{inverseinfactorization}
(t\ra u)^{-1}\la (t\la u)^{-1}=u^{-1},\quad (t\ra u)^{-1}\ra (t\la
u)^{-1}=t^{-1}
\end{eqnarray}
It  was shown in \cite{BGM} that $D(H)$ is a cocycle twist of the
double $D(X)=k(X)\lcross kX$, meaning in particular that its
category of modules is equivalent to that of $X$-crossed modules
in the sense of Whitehead.

Next, we need the notion of a bicovariant differential calculus
over any Hopf algebra $A$. A differential calculus over any
algebra $A$ is an $A-A$-bimodule $\Omega^1$ and a linear map
$d:A\to \Omega^1$ such that $\extd(ab)=a\extd b+(\extd a)b$ for
all $a,b\in A$ and such that the map $A\tens A\to \Omega^1$
defined by $a\extd b$ is surjective. In the Hopf algebra case we
require bicovariance in the sense that $\Omega^1$ is also an
$A-A$-bicomodule via bimodule maps and $\extd$ is a bicomodule
map\cite{W}, in which case one may identify $\Omega^1=A\tens
\Lambda^1$ where $\Lambda^1$ is the space of invariant 1-forms. It
forms a right $A$-crossed module (i.e. a compatible right
$A$-module and $A$-comodule or right module of the Drinfeld double
$D(A)$ in the finite dimensional case). The (co)action on
$\Omega^1$ from the left are via the (co)product of $A$, while
from the right it is the tensor product of that on $A$ and on
$\Lambda^1$. Then the classification amounts to that of
$\Lambda^1$ as quotient crossed modules of $A^+=\ker\eps\subset
A$.  Also, a calculus is irreducible (more precisely one should
say `coirreducible') if it has no proper quotients.  Then as in
\cite{MAj1} we actually classify the duals $L=\Lambda^1{}^*$,
which we call `quantum tangent spaces', as irreducible crossed
submodules of $H^+=\ker\eps\subset H$ under $D(H)$, where $H$ is a
Hopf algebra dual to $A$. Finally, we note that the category of
$A$-crossed modules is a braided one (since the Drinfeld double is
quasitriangular) and hence there is an induced braiding
$\Psi:\Lambda^1\tens\Lambda^1\to\Lambda^1\tens\Lambda^1$ which can
be used to define an entire `exterior algebra'
$\Omega(A)=A\tens\Lambda$. The invariant forms $\Lambda$ are
generated by $\Lambda^1$ with `antisymmetrization'
relations\cite{W} defined by $\Psi$.  We will use these notations
and concepts throughout the paper.

\section{Classification of differentials by conjugacy classes
in $X$}\label{classisection}

In this section we provide a concise but self-contained account of
the classification theory in \cite{BM}. We unfortunately need to
recall it in detail before we can derive the Cartan calculus
associated to each classification datum in Section~3. We will,
however, take the opportunity to reformulate the theory of
\cite{BM} more directly in terms of conjugacy classes and to
change to what are now more standard left-invariant forms. This is
not a routine left-right reversal of all formulae as the
bicrossproduct is not itself being reversed, and in fact leads to
cleaner results.

\subsection{Modules of the quantum double of a
bicrossproduct}\label{moduleofD(H)}

According to the Woronowicz theory \cite{W} in the form recalled
above, the first step to the classification is to understand the
$D(H)$-modules where $H=kM\rlbicross k(G)$, and in particular the
canonical one on $H^+$. We begin by recalling what is known about
these, from  \cite{BGM,BM} but with a necessary switch from left
to right modules. This is again not routine, but we omit the
proofs. Note that a $D(H)$ right module means a compatible right
module of $H$ and left module of $H^*$ (or right module of
$H^{*\rm op}$).

\begin{prop}\cite[Prop. 4.1]{BGM} \label{prop-D(H)-module}
The right modules of $D(kM\rlbicross k(G))$ are in one-one
correspondence with vector spaces $W$ which are:\\

(i)\quad $G$-graded right $M$-module such that $|w\ra t|=t^{-1}\la
|w|$, for  all $t \in M$, where $|\;\;|$ denotes the $G$-degree of
a homogeneous element $w \in W$.

(ii)\quad $M$-graded left $G$-module such that
 $\<u\la w\> =\<w\>\ra u^{-1}$, for  all
$u \in G$, where $\< \;\;\>$ denotes the $M$-degree of
 a homogeneous element $w \in W$.

(iii)\quad Bigraded by $G,M$ together and mutually ``cross
modules'' according to
\begin{eqnarray*}
\<w\ra t\>= t^{-1}\<w\>(t^{-1}\ra |w|)^{-1},\quad |u\la
w|=(\<w\>\la u^{-1})^{-1}|w|u^{-1}
\end{eqnarray*}

(iv)\quad $G-M$-''bimodules'' according to
\begin{eqnarray*}
\left((t^{-1} \ra |w|u)\la u^{-1}\right)\la (w\ra t)= (u^{-1}\la
w)\ra \left( t \ra (t^{-1}\<w\>\la u)\right)
\end{eqnarray*}
The corresponding action of the quantum double is given by
\begin{eqnarray*}
w\ra (t\otimes \delta _{v})=\delta_{t^{-1}\la |w|,v} w\ra t,\quad
(\delta _{s} \otimes u)\la w =\delta_{s,\<w\>\ra u^{-1}} u\la w
\end{eqnarray*}
and the induced braiding is
\begin{eqnarray*}
\Psi_{L,W}(l \otimes w)= w\ra (\<l\>^{-1}\ra |w|^{-1})^{-1}\otimes
(\<l\>^{-1}\la |w|^{-1})^{-1} \;\la l
\end{eqnarray*}
\end{prop}

In particular, $D(H)$ acts on $H$ by the standard right quantum
adjoint action of $H$ and by the left coregular action of $H^*$:
\[
g\ra h=Sh_{(1)}gh_{(2)},\quad a\la h = \sum h_{(1)}\otimes
<h_{(2)},a>\] where $g,h \in H$, $a \in H^*$, and $\Delta
h=h_{(1)}\tens h_{(2)}$ is the Sweedler notation. A routine
computation from the Hopf algebra structure of $kM\rlbicross k(G)$
yields these as
\[ (s\otimes \delta_u)\ra (t\otimes \delta_v)= \delta_{s\la
u,u(t\la v)^{-1}}(t''s{t''}^{-1} \otimes \delta_{t''\la u}),\quad
t''=(t\ra v)^{-1}\ra u^{-1}\]
\[
(\delta_t\otimes v) \la (s\otimes \delta_u) =   \delta_{t,s\ra
uv^{-1}}( s\otimes \delta_{uv^{-1}}).\] Comparing these with the
form of the actions in Proposition~\ref{prop-D(H)-module} we find
easily that the gradings, the $M-G$ actions for the right
canonical representation of $D(kM\rlbicross k(G))$ on
$W=kM\rlbicross k(G)$ and the induced braiding are \[ |s\otimes
\delta_u|=(s\la u)^{-1}u,\quad <s\otimes \delta_u>=s\ra u\] \[
(s\otimes \delta_u)\ra t =\bar {t}s{\bar {t}}^{-1} \otimes
\delta_{\bar {t}\la u}, \quad v\la (s\otimes \delta_u)=s\otimes
\delta_{uv^{-1}}\] \[ \Psi (s\otimes \delta_u \otimes t \otimes
\delta_v)= s''t{s''}^{-1} \otimes \delta_{s''\la v} \otimes s
\otimes \delta_{u(s''\la v)^{-1}(s'' t\la v)}\] where $ \bar {t}=
t^{-1}\ra (s\la u)^{-1}$ and $s''=(s\ra u)^{-1}\ra v^{-1}$.

Following the spirit of \cite{BGM} we can also give right
$D(H)$-modules in terms of the right modules of the quantum double
$D(X)=k(X)\lcross kX$ of the group $X$, where the action is by
$\Ad$. Explicitly, its Hopf algebra structure is \[ (\delta_x
\otimes y)(\delta_a \otimes b)=\delta_{y^{-1} xy,a}(\delta_x
\otimes yb),\quad \Delta (\delta_x \otimes
y)=\sum\limits_{ab=x}\delta_a \otimes y \otimes \delta_b \otimes
y\] and suitable formulae for the counit and antipode. It was
shown in \cite{BGM} that there is an algebra isomorphism $\Theta :
D(H)\longrightarrow D(X)$ defined by
\begin{eqnarray}\label{isomtheta}
\Theta (\delta_s \otimes u\otimes t\otimes
\delta_v)=\delta_{u^{-1}s^{-1}(t\la v)u} \otimes u^{-1}(t\ra v)
\end{eqnarray}
A straightforward computation shows that its inverse is
\begin{eqnarray*}
\Theta^{-1}(\delta_{su}\otimes tv)=\delta_{s^{-1}\ra (t\la
v)}\otimes
 (t\la v )^{-1}\otimes (t\ra \alpha )\otimes \delta_{\alpha ^{-1} v}
\end{eqnarray*}
where $\alpha =t^{-1} \la u^{-1}(s^{-1}t \la v)$. Hence $D(H)$ and
$D(X)$ modules correspond under these isomorphisms.

On the other hand, it is known that $D(X)$-modules $W$ are nothing
other than crossed modules in the sense of Whitehead, see
\cite{MAjFounda}, i.e. given by $X$-graded $X$-modules with
grading $||\ ||$ and (right) action $\tilde\ra$, say, compatible
in the sense
\begin{eqnarray}\label{gradingfacto}
 ||w \tilde{\ra} x||=x^{-1}||w||x
\end{eqnarray}
for all $x\in X$ acting on homogeneous $w\in W$. The corresponding
action is of course
\begin{eqnarray}\label{actionD(X)}
w \tilde{\ra}(\delta_x \otimes y) = \delta_{x,||w||} (w
\tilde{\ra} y),\quad \forall x,y\in X.
\end{eqnarray}
It is easy to see that the correspondence with the gradings and
actions in Proposition~\ref{prop-D(H)-module} is
\eqn{XMG-grad}{||w||=\<w\>^{-1}|w|,\quad w\tilde{\ra}
us=(u^{-1}\la w) \ra (s^{-1}\ra |u^{-1}\la w|^{-1})^{-1}} $\forall
w \in W,\  us \in X$.

Therefore the canonical representation of $D(H)$ that we are
interested in can be identified with such an $X$-crossed module.
Before giving it, following \cite{BM}, we identify the vector
space $kX$ spanned by $X$
 with the vector space $W=kM\rlbicross k(G)$ via
 $ vt \equiv t\otimes \delta_v$. Then

\begin{prop}\cite{BM}\label{Xcross-can} The right canonical
representation of $D(H)$ can be identified with $kX$ as an
$X$-crossed module

\begin{eqnarray}\label{actionD(X)onkX}
 {}\quad ||vt||=||t\otimes \delta_v||=v^{-1}t^{-1}v,\; vt \tilde{\ra}us
= (\tilde{s}\la vu)(\tilde{s}t\tilde{s}^{-1}),\;
\tilde{s}=s^{-1}\ra (vu)^{-1}.
\end{eqnarray}
\end{prop}

Finally, we are actually interested in the canonical action not on
$H$ but on $H^+$. This is\cite{MAj1} the right quantum adjoint
action as before and $h\ra a =\sum h _{(1)}<h_{(2)},a>-<a,h>1$ for
all $h\in H^+$. It is arranged so that the counit projection to
$H^+$ is an intertwiner. Therefore in our case \eqn{intertwinner}{
\Pi :kX \rightarrow H^+,\quad vt\mapsto t\otimes \delta_v-\eps(
t\otimes \delta_v)1=t\otimes \delta_v- \delta_{v,e}} is an
intertwiner between this action $\ra$ (viewed as an action of
$D(X)$) and the action (\ref{actionD(X)}) defined by the crossed
module structure.

\subsection{Quantum tangent spaces in $kM\rlbicross k(G)$}

We are now ready briefly to reformulate the classification
\cite{BM} for the quantum tangent spaces $L\subset H^+$ of
bicrossproduct quantum groups $H=kM\rlbicross k(G)$. The minor
technical innovation is to rework the theory in terms of a subset
$Z\subset X$ stable under conjugation in $X$. Here
\eqn{imageunder||.||}{Z=\image(\CN),\quad \CN:X\to X,\quad \CN
(vt)=||vt||} is manifestly stable since
$(us)^{-1}||vt||us=||vt\tilde{\ra} us||$, for all $vt,us \in X$ as
an expression of the $X$-crossed module structure of $kX$ in
Proposition~\ref{Xcross-can}. Working with $Z$ is obviously
equivalent to working as in \cite{BM} with the quotient $X/\sim$,
where $x\sim y$ if $\CN(x)=\CN(y)$. Moreover, orbits under
$\tilde\ra$ as in \cite{BM} now correspond to conjugacy classes in
$Z$. We denote respectively by $X_z$ and $C_z$ the centralizer and
the conjugacy class in $X$ of an element $z \in Z$. Clearly, $Z$
is the partition into conjugacy classes of its elements. All
results in this section are along the lines of \cite{BM} with such
differences.

\begin{prop}\label{j_zasX_zmodule}
For each $z\in Z$  we set $J_z =k\CN^{-1}(z)$.\\
(i)\quad The space  $J_z$ is a right $X_z$ representation\\
(ii) \quad $M_{C_z}=\bigoplus_{z' \in C_z}J_{z'}\subset  kX$ is a
subrepresentation under the right action of $k(X)\lcross kX$ from
Proposition~\ref{Xcross-can}. Moreover, $kX=\bigoplus_{ C_z}
M_{C_z}$ is the decomposition of $kX$ into such
subrepresentations.
\end{prop}

\proof Statement $(i)$ is immediate. We now prove $(ii)$. The
action of $\delta_z \in D(X)$ denoted by $\tra \delta_z$ is a
projection operator that projects $kX$ onto $J_z$. Then we have
$kX=\bigoplus_{z\in Z}J_z$. Since $Z$ is a partition by the
conjugacy classes $C_z$, we have
$$kX=\bigoplus_{C_z}\bigoplus_{z\in C_z}J_z= \bigoplus_{C_z}M_{C_z}.$$
For a chosen conjugacy class $C$, let us set $$\pi_C
=\sum\limits_{z \in C}(\tra \delta_z).$$  The operator $\pi _C$ is
a projection of $kX$ onto $M_C$. To show that $M_C$ is a right
$D(X)$ representation, it is enough to show that the action $\tra
(\delta_x \otimes y)$ of any $\delta_x \otimes y \in D(X)$
commutes with $\pi_C,$ i.e $\pi _C \circ (\tra (\delta_x\otimes
y))=(\tra (\delta_x\otimes  y))\circ \pi_C$. This  is an easy
computation using the crossed relation
$y\delta_z=\delta_{yzy^{-1}}y$ in $D(X)$. \eproof

From now we fix a conjugacy class $C_0$ of an element $z_0 \in Z$
, denote by $X_0$ the centralizer of $z_0$ in $X$ and set
 $J_0=k\CN^{-1}(z_0)$.

\begin{prop} \label{centralizermodule}
Let $J_0=J_1 \oplus J_2...\oplus J_n$ be the decomposition into
irreducibles under the action of $X_0.$ For each $z =\bar
z^{-1}z_0 \bar z \in C_0$ , we set $J_{i\bar z}=J_i \tra \bar z$ (
this does not depend on the choice of $\bar z$), then
$$M_i=\oplus_{z \in C_0}J_{i\bar z}\subset M_{C_0},\quad  1 \leq
i \leq n$$ are irreducible subrepresentations under the right
action of $k(X)\lcross kX$. Moreover, $M_{C_0}=\oplus_i M_i$ is a
decomposition of $M_{C_0}$ into irreducibles.
\end{prop}

\proof First of all we prove that $J_{i\bar z}$ does not depend on
the choice of $\bar z$. Indeed suppose that $\bar z^{-1}z_0 \bar z
=y^{-1}z_0 y=z'$. Then we have $y\bar z^{-1}\in X_0$ which implies
that $J_i \tra y\bar z^{-1}=J_i$ hence
\[J_{iy}=J_i \tra y= (J_i \tra y\bar z^{-1})\tra \bar z=J_i \tra
\bar z=J_{i\bar z} \] Next, by equivariance of $\CN$ one shows
easily that $J_{i\bar z}\cap J_{iy}=\{0\}$ if $\bar z^{-1}z_0 \bar
z \neq y^{-1}z_0 y$. So $M_i$ as shown is a direct sum.  Reasoning
as in \cite{BM} with suitable care, one shows that $M_i$ is a
right $k(X)\lcross kX$-module: the essential steps are the
following: Let $P_i: J_0 \rightarrow J_0$ be a right $X_0$-map
which projects to $ J_i \subset J_0$ with all other $J_j$
contained in its  kernel. Let us define the map $Q_i: M_{C_0}
\rightarrow  M_{C_0} $ by
\begin{eqnarray}
Q_i =\sum\limits_{z \in C_0}(\tra \bar z)\circ P_i \circ (\tra
\bar z^{-1})\circ (\tra \delta _{z})
\end{eqnarray}
It is clear that  $Q_i$ is a projection onto $M_i$. The similar
computations  as in \cite{BM} yield
\[Q_i \circ (\tra (\delta_a \otimes b))=(\tra (\delta_a \otimes b))\circ Q_i,\]
proving that $M_i$ is a $D(X)$-module. Moreover it is clear that
$\sum\limits_{z \in C_0}(J_0\tra \bar z)=\sum\limits_i M_i$ and
since $Q_iQ_j=0$ for $i\neq j$ we have $\sum\limits_{z \in
C_0}(J_0\tra \bar z)=\bigoplus_i M_i$. Finally one may verify that
$M_i$ is irreducible as $D(X)$-right module. \eproof

We therefore have a decomposition of $kX$ into irreducibles, for
every choice of conjugacy class $C$ of an element $z_0 \in Z$ and
every irreducible subrepresentation of the centralizer of $z_0$ in
$X$. The converse also holds:

\begin{prop}  \label{M_0fromM}
Let $\CM\subset kX$ be an irreducible right $k(X)\lcross kX$
representation under the action from Proposition~\ref{Xcross-can}.
Then as vector space, $\CM$ is of the form
$$\CM=\bigoplus _{z \in C} (\CM_0 \tra \bar z)$$
For some conjugacy class $C$ in $X$ of $z_0 \in Z$ and some
irreducible subrepresentation $\CM_0 \subset J_{z_0}$ of the
centralizer $X_0$ of $z_0$ in the group $X$.
\end{prop}
\proof We choose $z_0 \in Z$ such that $\CM_0:= \CM\tra
\delta_{z_0}$ is nonzero. Hence $\CM_0\subset J_{z_0}$. Moreover
$\CM_0$ is a $X_0$ subrepresentation of $J_{z_0}$. Indeed let
$m\tra \delta_{z_0} , m \in  \CM$ be an element of $\CM_0$ and $g
\in X_0$. We note that $m\tra g \in \CM$ since $\CM$ is a
$D(X)$-module. We note also that $g
\delta_{z_0}=\delta_{(g^{-1}z_0g)}g=\delta_{z_0}g$. Hence
\[(m\tra \delta_{z_0})\tra g=m\tra \delta_{z_0}g
=m\tra g \delta_{z_0} \in \CM_0 \] which shows that $\CM_0 $ is a
$X_0$-subrepresentation of $J_{z_0}$. Next if $J_1$ is an
irreducible subrepresentation of $\CM_0 $ under the action of
$X_0$, then by the preceding proposition $\bigoplus_{z \in C}(J_1
\tra \bar z)\subset \CM $ is an irreducible right representation
of $D(X)$. And since $\CM$ is irreducible we have
 $\CM=\bigoplus_{z \in C}(J_1 \tra \bar z)$. Finally note that $J_1$
is in fact $\CM_0 ,$ so that $\CM_0 $ is irreducible as
$X_0$-module, indeed by  Proposition \ref{centralizermodule}, two
distinct subrepresentations $J_1$ and $J_2$ of $X_0$ should give
distinct irreducible subrepresentations $ \sum\limits_{z \in
C}(J_1 \tra \bar z)\subset \CM$ and $ \sum\limits_{z \in C}(J_2
\tra \bar z)\subset \CM$. This is not possible as
 $\CM$ is irreducible. \eproof

Application of $\Pi:kX\to H^+$ from Section~2.1 then tells us that
we obtain subrepresentations of $H^+$ under the action of $D(X)$
by projecting via $\Pi$ the subrepresentations of $kX$. We can now
give the total description of the irreducible quantum tangent
spaces of $H$. Then cf. \cite{BM},
\begin{theo}\label{tangentspacestheo}
The irreducible quantum tangent spaces $L\subset H^+$
are all given by the following 2 cases:\\
(a)\quad For a conjugacy class $C\neq \{e\}$ of an element $z_0\in
Z$, for each irreducible right subrepresentation $\CM_0\subset
J_{z_0}$ of the centralizer of $z_0$, we have an irreducible right
$D(H)$-module
 $\CM=\bigoplus_{z \in C}(\CM_0 \tra \bar z)$ and an
isomorphic irreducible right subrepresentation $L=\Pi (\CM)\subset
H^+$.\\
(b)\quad  For $C=\{e\}$,  $J_e=kG$, $X_e=X$ and for any nontrivial
nonzero irreducible right subrepresentation $\CM_0\subset kG$ we
obtain an irreducible right $D(H)$-module $\CM=\bigoplus_{z
=e}(\CM_0 \tra \bar z)=\CM_0$ and the isomorphic
$D(H)$-subrepresentation $L=\Pi(\CM_0) \subset H^+$.
\end{theo}
\proof These steps are the same as in \cite{BM}. Briefly, if
$\CM=\bigoplus_{z \in C}(\CM_0 \tra \bar z)$ is an irreducible
representation of the unprojected action then by equivariance, the
map $\Pi :\CM \rightarrow L$ is a map of representations. It is
surjective. If it is 1-1 then the two representations are
isomorphic. The unique case where  $\Pi$ is not 1-1 is where $\bar
1:=\sum\limits_{u\in G}u\in \CM$ i.e the case where $\CM=k\bar 1$
and hence $\Pi(\CM)=\{0\}$, since $\CM$ is irreducible. This case
is the one excluded in the theorem. Conversely if $L$ is an
irreducible right subrepresentation of $H^+$ under $k(X)\lcross kX
$ then the inverse image $\Pi ^{-1}( L)\subset kX$ is also a
representation of $k(X)\lcross kX$ and it contains $k\bar 1$. If
$L\neq 0$ then $\Pi ^{-1} (L)$ contains at least one other
irreducible representation $\CM$ such that $k\bar 1\oplus \CM
\subset \Pi ^{-1} (L)$ then $\CM$ must be of the form described
above and by irreducibility of $L$, $\Pi (\CM)=L$. \eproof

We note that the  element $z_0$ is not strictly part of the
classification of the differential calculi. In fact an irreducible
bicovariant differential calculus is defined by a conjugacy class
$C$ and a irreducible $D(X)$-subrepresentation $\CM\subset kX$
such that $||\CM||=C$, where $||\CM||$ denotes the set of images
by $||.||$ of homogeneous elements of $\CM.$ It does not depends
on the chosen element in $C$. In the other words if
$$\CM=\bigoplus _{z \in C} (\CM_0 \tra \bar z)$$
with $\CM_0$ an irreducible subrepresentation of $J_{z_0}$ under
the action of the centralizer of $z_0$ then for any $z_1 \in C$ we
can write also $\CM$ as
 $$\CM=\bigoplus _{z' \in C} (\CM_1 \tra \bar{z'})$$
where $\CM_1$ an irreducible subrepresentation of $J_{z_1}$ under
the action of the centralizer of $z_1$. This follows from
Proposition~\ref{M_0fromM}. Indeed giving $\CM=\bigoplus _{z\in C}
(\CM_0 \tra \bar z)$, and $z_1 \in C$, we set $\CM_1 =\CM\tra
\delta_{z_1}$. This is nonzero since $\CM_0\tra \bar {z_1}\subset
\CM\tra \delta_{z_1}$. $\CM_1 = \CM\tra \delta_{z_1}\neq 0$
implies by the proof of Proposition ~\ref{M_0fromM} that $\CM_1 $
is an irreducible subrepresentation of $k\CN^{-1}(z_1)$ under the
action of the centralizer $G_{z_1}$ of $z_1$ and moreover
  $\CM=\bigoplus _{z \in C} (\CM_1 \tra \bar z)$.

With this characterization of the quantum tangent spaces in terms
of conjugacy classes  and centralizers, we recover the well known
cases where $H=kM$ or $H=k(G)$:

\begin{prop}\label{tensorcalc}
(i) \quad Set $G=\{e\}$ then $X=M$ and
Theorem~\ref{tangentspacestheo} recovers the usual classification
of the irreducible bicovariant calculi on
$A=k(M)$ by nontrivial conjugacy classes in $M$.  \\
(ii) \quad Set $M =\{e\}$ then $X=G$ and we recover the
classification for calculi on $A=kG$ by nontrivial irreducible
subrepresentations  $V\subset kG$ under the regular right action
of $G$ on itself as in \cite{MAj1}.
\\
(iii)\quad Set $X=G\times M$ with trivial actions. Then
$A=k(M)\tens kG$. An irreducible bicovariant calculus on $A$ is
defined by an conjugacy class $C$ in $M$ and an irreducible
subrepresentation $V\subset kG$ under the regular right action of
$G$ on itself, with at least one of $V,C$ nontrivial. Here
$\CM=V.C$.
\end{prop}
\proof For case $(i)$ the action of $X$ on $kX$ in
Proposition~\ref{Xcross-can} is $t\tra s=s^{-1}ts$, $Z=M$. For any
conjugacy class $C_0$ of an element $t\in M$ we denote by $C$ the
conjugacy class of $t^{-1}$ and we have $J_t=k\{t^{-1}\},$ since
$||b||=b^{-1},\; \forall b\in M$. Hence the corresponding
irreducible subrepresentation  $\CM \subset kM$ under the action
of $D(H)$  is $\CM=\sum\limits_{z \in C}(J_t \tra \bar
z)=\sum\limits_{z \in C}k({\bar z}^{-1}t^{-1}\bar z)=kC $ and
$L=k\{a-e,a\in C\}$, i.e the basis of $L^*$ is labelled by a
conjugacy class as usual.

For case $(ii)$ the action of $X$ on $kX$ is $v\tra u=vu$,
$||v||=e, \;\forall v\in G$. Hence $Z=\{e\}$ so that we are in
case (b) of the theorem. Therefore the quantum tangent spaces
$L\subset H^+$ are isomorphic to the irreducible
subrepresentations  $V\subset kG$ as stated.

For case $(iii)$ we have $Z=M$. The action of $X$ on itself is
$vt\tra us=vu.s^{-1}ts.$ Let us consider a conjugacy class
$C_{t_0^{-1}}$ of $t_0^{-1}$ in $Z$. The centralizer of $t_0$ in
$X$ is $X_0=G.cent_M(t_0)$, where $cent_M(t_0)$ is the centralizer
of $t_0$ in $M$, $J_0=k\CN^{-1}(t_0^{-1})=kG.t_0.$  The action of
$X_0$ on $J_0$ is \[v.t_0\tra us =vu.t_0,\quad \forall v\in G,
us\in X_0\] which leads to $\CM_0$ of the form $\CM_0=V.t_0$,
where $V$ is as mentioned, hence
$$\CM=\bigoplus_{t\in C_{t_0^{-1}}} (V.t_0\tra \bar t)=V.C_{t_0}$$
where $C_{t_0}$ is the conjugacy class of $t_0$ in $M.$ \eproof

The calculus in case $(iii)$ is a product of calculi on $G,M$ for
the cases $(i)$ and $(ii)$ and has the product of their
dimensions.

\section{Cartan calculus on $k(M)\lrbicross kG$}\label{cartancalculus}

We are now ready to proceed to our main results. Let
$A=k(M)\lrbicross kG$ be the dual of $H=kM\rlbicross k(G).$ Our
goal is to find an explicit description for the calculus
corresponding to each choice of classification datum. This amounts
to a description of the differential forms and the commutation
relations with functions and  $\extd$, i.e. a `Cartan calculus'
for the associated noncommutative differential geometry.

We fix a conjugacy class $C$ of an element $z_0 \in Z$, an
irreducible right subrepresentation $\CM_0\subset J_{z_0}$ of the
centralizer of $z_0$, and  the corresponding  nontrivial
irreducible right $D(H)$-module $\CM=\bigoplus_{z \in C}(\CM_0
\tra \bar z)$  as in Theorem \ref{tangentspacestheo} above. For
each $z\in C$ we fix one element $\bar z$ so that $z=\bar
z^{-1}z_0\bar z$ and we set $\bar C =\{\bar z|\ z \in C\}$. As we
saw above, $\CM=\bigoplus_{\bar z \in \bar C}(\CM_0 \tra \bar z)$.
We now choose a basis $(f_i)_{i\in I}$ of $\CM_0$ ($I$ is finite)
and set
\begin{eqnarray*}
f_{iz}:=f_i \tra \bar z
\end{eqnarray*}
We recall that here $\tra$ is the action of $X$ on itself defined
in Proposition~\ref{Xcross-can}.

\begin{lemma}\label{basesofL}
The vectors $(f_{iz}),i\in I$  form a basis of  $\CM$ with
homogeneous $X$-degree $z$.
\end{lemma}
\proof: By definition  it is clear that $(f_{iz})$  generate $\CM$
since $(f_i)$ generate $\CM_0 $. Using the direct sum in the
decomposition of $\CM$ and the fact that $(f_i)$ are linearly
independent, one checks easily that $(f_{iz})$ are linearly
independent too. By definition, $(f_i)$ are homogeneous of degree
$z_0$. This implies that each $f_{iz}$ is homogeneous of degree
$z$ since for homogeneous $w$, $||w\tra x||=x^{-1}||w||x,\forall
x\in X.$ \eproof

In what follows, we  identify $\CM$ with the quantum tangent space
$L$ as isomorphic vector spaces via  $\Pi$. The dual $\Lambda^1$
of $L$ is equipped with the dual basis $(e_{iz})$ of the basis
$(f_{iz})$.

To simplify we relabel these basis  by $(e_a)_{a\in I}$ and
$(f_a)_{a\in I}$ respectively for the space of invariant 1-forms
and the quantum tangent space. We recall the factorization
(\ref{XMG-grad}) of an $X$-grading into an $G$-grading $|\ |$ and
an $M$ grading $\<\ \>$.
\\

We are now ready to follow the Woronowicz construction explained
in the preliminaries to build $(\Omega^1,\extd)$ as a differential
bimodule, namely we set $\Omega^{1}(A)=A\otimes \Lambda^1$,
\begin{eqnarray}\label{dextformula}
\quad \extd a=\sum (\id \otimes \Pi _{\Lambda^1})(a_{(1)}\otimes
(a_{(2)}-\eps (a_{(2)}))),
\end{eqnarray}
\begin{eqnarray}\label{commurelationsformula}
a.x=a\otimes x,
 \quad x.a =\sum a_{(1)}\otimes x\ra a_{(2)}
\end{eqnarray}
for all $a\in A$ and $x\in\Lambda^1$, where $\Pi_{\Lambda^1}$
denotes the projection of $A^+$ on $\Lambda^1$ adjoint to the
injection $L\subset H^+$. We have the following :

\begin{theo}\label{cartanrelations}
With the chosen basis of $L$ as above, the differential calculus
in Theorem~\ref{tangentspacestheo} is explicitly defined by:\\

(i)\quad The left $A$-module of 1-forms $\Omega^{1}(A)=A\otimes
\Lambda^1$.\\

(ii)\quad The right module structure according to commutation
relations between ``functions'' and 1-forms:
\[
e_a\delta_s=\delta_{s\<f_a\>^{-1}}e_a, \quad e_au=(\<f_a\>\la u)
e_a*u\] where
\[e_a*x=\sum\limits_{b\in I}<e_a,f_b\tra x^{-1}>e_b, \quad \forall x\in X.\]
is the right action of $X$ on $\Lambda^1$
adjoint to the left action $x*f_a:= f_a\tra x^{-1}$ on $L$.\\

$(iii)$ \quad The exterior differential:
\begin{eqnarray*}
\extd\delta_s= \sum\limits_{a} <\delta_{\<f_a\>},f_a>
(\delta_{s\<f_a\>^{-1}}-\delta_s)e_a
\end{eqnarray*}
\begin{eqnarray*}
\extd u=\sum\limits_{a}<\delta_{\<f_a\>},f_a> (\<f_a\>\la
u)e_a*u-\sum\limits_{a} <\delta_{\<f_a\>},f_a> ue_a
\end{eqnarray*}
where  $<\delta_{vt},f_a>$ for all $vt\in X,$ is the pairing
between $k(X)$ and its dual $kX$.
\end{theo}

\proof We first of all note the following facts easily obtained
from (\ref{gradingfacto}) and the factorization of $X$ grading in
(\ref{XMG-grad}) and that we freely use in the proof:
\[ \<f_a\tra u\>=\< f_a\>\ra u,\quad \forall u\in G,\quad
 |f_a\tra s|=(s^{-1}\la |f_a|^{-1})^{-1},\quad \forall s\in M.\]

We note also that the right action of $A$ in
(\ref{commurelationsformula}) is the restriction of the action of
$D(A^*)$ on $\Lambda^1$, we view it via the isomorphism $\Theta$
as action of $D(X)$ on $\Lambda^1$, adjoint to a left action of
$\Theta( A)$ on $L\cong \CM \subset kX.$ Clearly equation
(\ref{actionD(X)onkX}) expresses both right action of
$\Theta(A^*)$ and left action of $\Theta(A)$ on $\CM$, thus:
\begin{eqnarray}\label{actofAonlambda1formula}
e_a\ra \delta_b \tens u=:e_a \tra \Theta(\delta_b \tens u)
=\sum\limits_{c\in I} <e_a,f_c\tra\Theta( \delta_b \tens u)>e_c
\end{eqnarray}
for all $\delta_b \tens u \in A.$ On the other hand
\[\Theta(\delta_s \tens u)=\sum\limits_{v\in G}\delta_{u^{-1}svu}\tens
u^{-1}\] then by (\ref{relationsofAa}) and
(\ref{commurelationsformula}) we have
\[e_a\delta_s=\sum\limits_{v\in G, b\in M} \delta_{sb^{-1}}\tens
(e_a\tra \delta_{b^{-1}v}).\] We compute
\[e_a\tra \delta_{b^{-1}v}=
\sum\limits_{c\in I} <e_a,f_c\tra  \delta_{b^{-1}v}>e_c
=\sum\limits_{c\in
  I}\delta_{||f_c||,b^{-1}v}\delta_{c,a}e_c=\delta_{\<f_a\>,b}\delta_{|f_a|,v}e_a\]
so that
\begin{eqnarray*}
e_a\delta_s= \sum\limits_{v,b} \delta_{sb^{-1}}\tens
\delta_{\<f_a\>,b}\delta_{|f_a|,v}e_a= \delta_{s\<f_a\>^{-1}}\tens
e_a
\end{eqnarray*}
Next, let $u\in G\subset A.$ From (\ref{relationsofAa}) and
(\ref{commurelationsformula}) again, we have
\[e_au=\sum\limits_{b\in M}(b\la u)\tens e_a \tra \Theta ( \delta_b
\tens u)\] To compute $e_a \tra \Theta ( \delta_b\tens u),$ we
first note that if we change the basis $(f_a)$ to $f'_a=:f_a\tra
u$ then its dual $(e_a)$ transforms as
$$e'_a=:e_a*u=\sum\limits_{c\in I} <e_a, f_c\tra u^{-1}>e_c.$$
Then
\begin{eqnarray}\label{explicitactofAonLambda1}
e_a \tra \Theta ( \delta_b\tens u) &=& \sum\limits_{c\in I, v\in
G}<e_a, f'_c\tra \delta_{u^{-1}b^{-1}vu}\tens
u^{-1}>e'_c \nonumber \\
&=& \sum\limits_{c\in I, v\in G} \delta_{||f'_c||,u^{-1}b^{-1}vu}
<e_a,f'_c \tra u^{-1}>e'_c \nonumber\\
&=& \sum\limits_{c\in I, v\in G} \delta_{\<f'_c\>,b\ra
  u}\delta_{|f'_c|,
u^{-1}v^{-1}(b\la u)}
<e_a,f_c>e_c*u \nonumber\\
&=& \delta_{\<f_a\>,b}e_a*u
\end{eqnarray}
from which we deduce $e_au=(\<f_a\>\la u)e_a*u$ as required.

We now prove the formulae for differentials. Writing $\bar
a=a-\eps(a)\in A^+,$ we write the projection as
\[\Pi_{\Lambda ^1}(\bar a)=
\sum\limits_{c\in I} <\Pi_{\Lambda ^1}(\bar a), f_c>e_c=
\sum\limits_{c\in I} <\bar a, i(f_c)>e_c=\sum\limits_{c\in I}
<\bar a, \Pi(f_c)>e_c\] where $i$ is the injection $L\subset H^+$
which in our case, viewing $L$ as $\CM\subset kX$, is just the
restriction on $\CM$ of the map $\Pi: kX\to H^+$ in
(\ref{intertwinner}). Denoting its adjoint map $\Pi^*: H^* \to
k(X)$ we have therefore $\Pi_{\Lambda^1}(\bar a)=\sum\limits_{c\in
I} <\Pi^*(\bar a), f_c>e_c$. In our case,
\begin{equation}\label{pi*}
\Pi^*(\delta_{s}\tens u-\delta_{s,e}1\tens e)=
\delta_{us}-\delta_{s,e}\sum\limits_{t\in M}\delta_{e.t}
\end{equation}
Let $s\in M$. From (\ref{dextformula}) we have
\begin{eqnarray*}
d\delta_s&=&\sum\limits_{b\in M}\delta_{sb^{-1}}\tens
\Pi_{\Lambda^1}(\delta_b-\delta_{b,e}1_{A})=\sum\limits_{b\in
M}\delta_{sb^{-1}}\tens \delta_{e.b}-\delta_{b,e}\sum\limits_{t\in
M}\delta_{e.t} \\
&=&\sum\limits_{b\in M, c\in I}\delta_{sb^{-1}}\tens <
\delta_{e.b}-\delta_{b,e}\sum\limits_{t\in
M}\delta_{e.t},f_c>e_c\\
&=&\sum\limits_{b\in M, a\in I}\delta_{sb^{-1}}\tens
(<\delta_{\<f_a\>},f_a>\delta_{b,\<f_a\>}-\delta_{b,e}<\delta_{\<f_a\>},f_a>)e_a\\
&=&\sum\limits_{ a\in
I}<\delta_{\<f_a\>},f_a>(\delta_{s\<f_a\>^{-1}} -\delta_s)e_a,
\end{eqnarray*}
where we used
\eqn{deltabfa}{<\delta_{e.b},f_a>=<\delta_{\<f_a\>},f_a>\delta_{b,\<f_a\>},\quad
\forall b\in M,} as one may see by expanding
$f_a=\sum\alpha_a^{vt}vt$, say. This pairing also equals
$<\delta_{ub},f_a\tilde\ra u>$ for all $u\in G$ since
$<vt\tilde\ra
u,\delta_{ub}>=<vt,\delta_{e.b}>=\delta_{v,e}\delta_{t,b}$. Hence,
using (\ref{relationsofAa}) and (\ref{dextformula}), we similarly
have
 \begin{eqnarray*}
du&=&\sum\limits_{b\in M}(b\la u)\tens
\Pi_{\Lambda^1}(\delta_b\tens
u-\delta_{b,e}1_{A})\\
&=&\sum\limits_{b\in M, a\in I}(b\la u)\tens< \Pi^*(\delta_b\tens
u-\delta_{b,e}1_{A}),f'_a>e'_a\\
&=&\sum\limits_{a\in I, b\in M}(b\la u)<\delta_{ub},f'_a>e'_a
-\sum\limits_{a\in I, t\in M}u<\delta_{e.t},f_a>e_a\\
&=&\sum\limits_{a\in I}<\delta_{\<f_a\>},f_a>(\<f_a\>\la u)e_a*u
-\sum\limits_{a\in I}<\delta_{\<f_a\>},f_a>ue_a.
\end{eqnarray*}
 This ends the proof of Theorem
\ref{cartanrelations}. \eproof

\begin{corol}\label{theta}
All irreducible bicovariant differential calculi on a
bicrossproduct $A=k(M)\lrbicross kG$ are inner in the sense
\begin{eqnarray*}
 da=[\theta,a],\quad \forall a \in A,\quad{\rm where}\quad
 \theta=\sum\limits_{a\in I}c_ae_a,\quad
 c_a=<\delta_{\<f_a\>},f_a>.
\end{eqnarray*}
\end{corol}
\proof The relations   $\theta \delta_s -\delta_s
\theta=d\delta_s$ and  $\theta u-u\theta=du$ are obtained  from
the definitions in Theorem \ref{cartanrelations}. \eproof

Once the first order differential calculus is defined explicitly,
we need also the braiding $\Psi$ induced on $\Lambda^1\otimes
\Lambda^1$ (then on $\Omega^1(A)\otimes_A \Omega^1(A)$) to
determine $\Omega^n(A), n\geq 2$. Thus $\Omega^2(A)=A\tens
\Lambda^2$ where $\Lambda^2$ is the space of invariant 2-forms
defined as the quotient of $\Lambda^1\tens \Lambda^1$ by
$\ker(\id-\Psi)$.

\begin{prop}\label{braidingonVXV}
The braiding $\Psi$ induced on $\Omega^1(A)$ by the action of the
quantum double $D(A^*)$ is given by
\begin{eqnarray*}
\Psi(e_a\otimes e_b)&=&e_b*(\<f_{a}\>\ra
 |f_b|)^{-1}\otimes e_a*|f_b|.
\end{eqnarray*}
\end{prop}
\proof The formula of the braiding on a basis $(e_a\otimes e_b)$
of $\Lambda^1\otimes \Lambda^1$ is
 \begin{eqnarray}\label{braidingformula}
\Psi(e_a\otimes e_b)=\sum\limits_i \beta^i\la e_b \otimes e_a \ra
\alpha_i
\end{eqnarray}
where $(\alpha_i)$ is a basis of $A$ with dual basis $(\beta^i)$
of $A^*$. The right action of $A$ on $\Lambda^1$ is given by
(\ref{explicitactofAonLambda1}) i.e.
\begin{eqnarray}\label{actionofAonV}
e_a\ra (\delta_t \otimes v)=\delta_{t,\<f_{a}\>}e_a*v
\end{eqnarray}
for all $t\in M$ and $u\in G$. We now compute $(t\otimes \delta_v)
\la e_a$, the adjoint of the action  $f_a\ra(t\otimes \delta_v)$.
 We have
\begin{eqnarray*}
\Theta(t\otimes \delta_v)=\sum\limits_{s\in M}\delta_{s^{-1}(t\la
v)}\otimes t\ra v
\end{eqnarray*}
so that
\begin{eqnarray*}
f_b\ra (t\otimes \delta_v)&=&\sum\limits_{s\in M}f_b \tra
(\delta_{s^{-1}(t\la v)}\otimes (t\ra v)) = \sum\limits_{s\in
M}\delta_{||f_b||,s^{-1}(t\la v)}f_b \tra (t\ra v).
\end{eqnarray*}
Then considering the basis $f''_b=f_b\tra (t\ra v)^{-1}$ whose
dual basis is $e''_b=e_b*(t\ra v)^{-1},$ we compute
\begin{eqnarray*}
(t\otimes \delta_v)\la e_a&=&\sum\limits_{b\in I}<e_a,f''_b
\tra (t\otimes \delta_v)>e''_b\\
&=&\sum\limits_{b\in I,s\in M}
\delta_{||f''_b||,s^{-1}(t\la v)}<e_a,f''_b\tra (t\ra v)>e''_b\\
&=& \sum\limits_{s\in M}\delta_{||f''_a||,s^{-1}(t\la v)}e_a*(t\ra
v)^{-1}.
\end{eqnarray*}
It is easy to check that
\[||f''_a||=(\<f_a\>\ra (t\ra v)^{-1})^{-1}((t\ra v)\la |f_a|^{-1})^{-1}\]
so that
\begin{eqnarray}\label{explicitleftactionofAonLamda1}
(t\otimes \delta_v)\la e_a =\delta_{(t\la v)^{-1},(t\ra v)\la
  |f_a|^{-1}}e_a*(t\ra v)^{-1}=\delta_{v,|f_a|}e_a*(t\ra v)^{-1}.
\end{eqnarray}
Finally, combining equations (\ref{actionofAonV}) and
(\ref{explicitleftactionofAonLamda1}) gives the formula for the
braiding as stated. \eproof

\begin{corol}\label{thetainv} The left-invariant 1-form $\theta$ obeys
$\theta\wedge \theta =0$ and is closed and nontrivial in the first
noncommutative de Rham cohomology $H^1$. The basic 1-forms obey
the Maurer-Cartan relations $d e_a=\{\theta,e_a\}$.
\end{corol}
\proof We need only to prove that $\theta$ is right-invariant (the
rest then follows by general arguments). This is equivalent to
invariance under the left action
(\ref{explicitleftactionofAonLamda1}) of $H$, which is a modest
computation. Alternatively, the relevant coaction on $\Lambda^1$
is the projection of the adjoint one on $A^+$. At least for
$\CC^X\ne\{e\}$, we have $\theta=-\Pi_{\Lambda^1}(\delta_e\tens
e-\sum_{t\in M}\delta_t\tens e)$ by (\ref{pi*})-(\ref{deltabfa})
and similar computations as there. This representative element of
$A^+$ is then more obviously $\Ad$-invariant. Since this coaction
also enters into $\Psi$, invariance then implies that
$\Psi(x\tens\theta)=\theta\tens x$ for any $x\in \Lambda^1$ and
hence that $\Psi(\theta\tens\theta)=\theta\tens\theta$. This is in
any case true when $\C^X=\{e\}$ since $\Psi$ is then the usual
flip. Hence $\theta\wedge \theta=0$ in the exterior algebra. On
the other hand, for the Woronowicz construction for any Hopf
algebra one may show that if the first order calculus is inner by
a left-invariant 1-form $\theta$ obeying $\theta\wedge \theta=0$
then the entire exterior calculus is inner, i.e. $d
\rho=[\theta,\omega\}$ for any form $\omega\in \Omega$. The graded
commutator here denotes commutator in degree 0 and anticommutator
in degree 1. Hence the last part of the Corollary is automatic. It
implies then that $d \theta=0$.

It remains only to show that $\theta$ is not exact. This is
actually true for any left-invariant 1-form on a left-covariant
calculus when the Hopf algebra is semisimple. Precisely such Hopf
algebras have a (say) right-invariant integral $\int:A\to k$ such
that $\int 1=1$ (for our bicrossproduct $A$ it is $\int
(\delta_s\tens u)=|M|^{-1}\delta_{u,e}$ as in \cite{MAjFounda}).
In this case suppose $d a\in \Lambda^1$ for some $a\in A$, so that
$\Delta_L(d a)=a_{(1)}\tens \extd a_{(2)}=1\tens \extd a$, then
$(\int a_{(1)})  \extd a_{(2)}=\int(1) \extd a=\int(a) d(1)=0$ by
right-invariance of the integral and $d(1)=0$. Hence $\theta$ is
necessarily nontrivial in the noncommutative de Rham cohomology.
 \eproof

\begin{prop}\label{groupcase} We recover the results known in
the cases of the group algebra and functions algebra of a finite
group.
\end{prop}\proof \quad (i)\quad  For $G=\{e\}$, $A=k(M)$ and
$L=\Pi(kC)$ for a conjugacy class $C$. Here a basis of $kC$ is
$(f_a=a)_{a\in C}$ since the action is $t\tra s=s^{-1}ts$.
Moreover $\<f_a\>=a$, $|f_a|=e$ and
$e_b*a^{-1}=e_{aba^{-1}},\;\;\forall a,b \in C.$ Then Theorem
\ref{cartanrelations}, Corollary \ref{theta} and Proposition
\ref{braidingonVXV} read \[\Omega^1(A)=A\otimes (kC)^*;\quad
e_a\delta_s=\delta_{sa^{-1}}e_a\]
\[ d(\delta_s)=\sum\limits_{a\in
C}<\delta_a,f_a>(\delta_{sa^{-1}}-\delta_s)e_a= \sum\limits_{a\in
C}(\delta_{sa^{-1}}-\delta_s)e_a\] which is exactly (\ref{dk(G)}
on a general function $f\in k(M)$. Moreover,
\[  \theta= \sum\limits_{a\in
C}<\delta_a,f_a>e_a=\sum\limits_{a\in C}e_a,\quad \Psi(e_a\otimes
e_b)=e_{aba^{-1}}\otimes e_a\] for $a,b\in C,$ $s\in M$ and
$R_a(f)(x)=f(xa)$, for all $x\in M$. (ii)\quad For $M=\{e\}$,
$A=kG$. $L$ is an irreducible subrepresentation of $kG$ under the
right multiplication in $G$. Let $(f_j)$ be a basis of $L$ with
the dual basis $(e_j)$. Here we have $||f_j||=e, \;\; \forall j.$
Then Theorem \ref{cartanrelations}, Corollary \ref{theta}  and
Proposition \ref{braidingonVXV} read
\[
\Omega^1(A)=A\otimes L^*;\quad e_iu =ue_i*u\]
\[du=\sum\limits_i<\delta_e,f_i>(ue_i*u-ue_i)=u\theta *(u-1)\]
\[
\theta=\sum\limits_i<\delta_e,f_i>e_i,\quad \Psi(e_i\otimes
e_j)=e_j\otimes e_i\] as in (\ref{dkG}). \eproof

Hence the Cartan calculus in theorem~\ref{cartanrelations} indeed
generalizes the ones on  the group algebra and on the algebra of
functions of a finite group. \eproof

\section{Differential calculi on cross coproducts $k(M)\rcrossco kG$.}

Now that we have the  Cartan calculus  for general bicrossproduct
Hopf algebras,  we specialize to the semidirect case where
$X=G\lcross M$ or $A=k(M)\rcrossco kG$, a cross coproduct. These
are the `coordinate' algebras of semidirect product quantum groups
$H$. In this case some further simplifications are possible.

We start with a general observation about the structure of $Z$ for
general $X=GM$. As usual, $u,v,g...$ are elements of $G$ and
$s,t,\tilde s...$ are those of $M$.

\begin{prop} (i) For general $X=G.M,$
the set $Z$ is given in terms of conjugacy classes $\CC^{M}$ of
$M$ by
$$Z=\bigcup_{\CC^{M}} \bigcup_{u\in G}(u^{-1}\CC^{M}u)$$
and for any fixed conjugacy class $\CC^{M} \subset M,$ the set
$$C^X=\bigcup_{u\in G}(u^{-1}\CC^{M}u)$$
 is a conjugacy class in $X$ .\\
(ii) In the semidirect case $X=G\lcross M$, the map
$$C^Z: \CC^{M}\longrightarrow  \bigcup_{u\in G}(u^{-1}\CC^{M}u)$$
from the set of conjugacy classes of $M$ to that of conjugacy
classes of $X$ contained in $Z$ is  one to one.
\end{prop}
\proof We first note that the map $C^Z$ is not one to one in
general (e.g.  for the $\Z_6.\Z_6$ example, $\CC_t^{M}$ and
$\CC_{t^{-1}}^{M}$ are different and have the same image through
$C^Z$).

In the semidirect case $X=G\lcross M$ this map is one to one since
\begin{eqnarray*}
(u^{-1}t_1u)=(v^{-1}t_2v)
&\Longleftrightarrow &(u^{-1}(t_1\la u)).t_1=(v^{-1}(t_2\la v)).t_2\\
&\Longrightarrow & t_1=t_2.
\end{eqnarray*}
 The other assertions  are easily obtained too. \eproof

\subsection{Canonical calculus for the case $X=G\lcross M$}\label{Gactstrivaily}

Now we specialize to the semidirect case $X=G\lcross M$. As we saw
above, an irreducible differential calculus on $A$ is defined by a
conjugacy class (of $t_0^{-1}\in M$ say) $\CC_0^{X}\subset Z$ and
a choice of an irreducible subrepresentation $\CM_0$ of
$J_0=k\CN^{-1}(t_0^{-1})$ under the action of the centralizer of
$t_0$ in $X$. In the semidirect case, we have

\begin{prop}\label{semidirectcase} {\ }

(i)\quad $\CC_0^{X} =\{u^{-1}(t\la u).t\in X,\quad  t\in
\CC_0^{M}, u\in G \}$

(ii)\quad$\CN^{-1}(t_0^{-1})=N_0.t_0$

\noindent where $N_0= \{u\in G,  t_0\la u=u\}$ is a subgroup of
$G$.

(iii)\quad The centralizer $X_0$ of $t_0$ in $X$ is
$X_0=N_0.cent(t_0)$.

(iv)\quad The action of $X_0$ on $J_0=kN_0.t_0$ is given by
\[v.t_0\tra us=(s^{-1}\la vu).t_0,\quad \forall v,u\in N_0, s\in
cent(t_0).\]

(v)\quad There is a canonical choice of $\CM_0$ (hence a canonical
choice of an irreducible calculus on $A$)  defined by a conjugacy
class $\CC_0^{M}\subset M$.
\end{prop}
\proof
 $(i)$, $(ii)$ and $(iii)$  are immediately obtained from the definition
of $\CC_0^{X},\CN^{-1}$ and  $G_0$. For $(iv)$, the action of
$G_0$ on $J_0=kN_0.t_0$ is given by
\[v.t_0\tra us=(s^{-1}\la vu).(s^{-1}t_0s)=(s^{-1}\la vu).t_0,
\quad v,u\in N_0, s\in cent(t_0).\] For $(v)$, the element
 \[m_0=\sum\limits_{v\in N_0} v.t_0\]
 generates a one-dimensional
trivial  $N_0.cent(t_0)$-irreducible subrepresentation of
$J_0=kN_0.t_0$, since for all $u.s \in N_0.cent(t_0)$
\begin{eqnarray*}
m_0\tra u.s&=&\sum\limits_{v\in N_0}v.t_0\tra
u.s=\sum\limits_{v\in N_0}(s^{-1}\la vu).t_0=\sum\limits_{v\in
N_0}v.t_0=0\end{eqnarray*} where the penultimate equality is by
freeness of $\la\circ R_u$,with $R_u$=right multiplication. Hence
if $t_0\neq e$ then
$$\CM=\bigoplus_{ z\in {C_0^X}}k.m_0\tra \bar z$$
 is the corresponding quantum tangent space with dimension $|C_0^X|$.
Hence, to any conjugacy class of $M$ (or any irreducible
differential calculus on $k(M))$ corresponds a canonical
irreducible differential calculus on $A$. Here, we made a
convention that the null calculus corresponds to  $t_0=e$. \eproof

As an important subcase, we consider now $X=G\lcross G$ where the
action is by conjugation. In this case $A=k(G)\rcrossco kG=D^*(G)$
is the dual of the quantum double of the group algebra $kG$. Then
Proposition \ref{semidirectcase} reads
\begin{corol}\label{corolG=G} When $M=G$ and $X=G\lcross G$
by conjugation, we have

(i) \quad $\CC_{0}^{X}=\bigcup_{s\in
\CC_0^{M}}(\CC_0^{M}s^{-1}).s$

(ii) \quad $\CN^{-1}(t_0^{-1})=cent(t_0).t_0$

(iii) \quad $X_0=cent(t_0).cent(t_0)$

(iv)\quad The action of $u.s \in cent(t_0).cent(t_0)$ on $v.t_0
\in \CN^{-1}(t_0^{-1})$ is
\[
v.t_0\tra u.s=s^{-1}vus.t_0=Ad_{s^{-1}}\circ R_u(v).t_0\]
\end{corol}
\proof We check easily that $N_0$ becomes $cent(t_0)$ and the
results stated  follow immediately from  Proposition
\ref{semidirectcase}. \eproof

Moreover, part $(v)$ of Proposition \ref{semidirectcase} says that
 any irreducible differential calculus on $k(G)$ extends to a
canonical irreducible differential calculus on $A=D^*(G)$. We
describe it explicitly. Let
\[C_0^M=\{s_0=t_0^{-1},s_1,...,s_N\}\]
be a conjugacy class (of $t_0^{-1}$) in $M=G$ and $C_0^X$ be the
corresponding conjugacy class of $t_0^{-1}$ in $X $ as above. For
each $0\leq i\leq N$, we fix $\bar {s_i}$ in $M$ such that
\[ s_i={\bar {s_i}}^{-1}t_0^{-1}\bar {s_i}.\]
To avoid confusion we use here the following notation: $s_i$ is
always in $M$ and we let $\und{\ }$ denote the identity map from
$M$ to $G$, so $\und {s_i}$ denotes the same element in $G$. As
usual, in any expression $g.t\in X$, we have $g\in G$ and $t\in
M.$ Then by $(i)$ of Corollary \ref{corolG=G}, each element
$z_{ij}$ of $C_0^X$ is of the form
\[z_{ij}=\und {s_i}\und {s_j}^{-1}.s_j, \quad s_i,s_j\in C_0^M\]
The elements $\bar {s_i}$ define $\overline{z_{ij}}\in X$ such
that $z_{ij}=\overline{z_{ij}}^{-1}t_0^{-1}\overline{z_{ij}}$ and
we have
\[\overline{z_{ij}}=\bar {s_i}\bar {s_j}^{-1}.\bar {s_j}.\]
Indeed if we set $g_{ij}=\bar {s_i}\bar {s_j}^{-1}.e=\und {\bar
{s_i}} \und {\bar {s_j}}$ then we have
\begin{eqnarray*}
\overline{z_{ij}}^{-1}t_0^{-1}\overline{z_{ij}}&=&\bar
{s_j}^{-1}g_{ij}^{-1}
t_0^{-1}(g_{ij}.\bar {s_j}) \\
&=&(\bar {s_j}^{-1}\la g_{ij}^{-1})(e.\bar {s_j}^{-1}t_0^{-1})
(g_{ij}.e)(e.\bar {s_j})\\
&=&(\bar {s_j}^{-1}\la g_{ij}^{-1})(\bar {s_j}^{-1}t_0^{-1}\la
g_{ij})
.\bar {s_j}^{-1}t_0^{-1}\bar {s_j}\\
&=&(\bar {s_j}^{-1}\la g_{ij}^{-1})(t_0^{-1}\la  g_{ij}).s_j\\
&=&\und {s_i}\und {s_j}^{-1}.{s_j}=z_{ij}
\end{eqnarray*}
We are now in position to compute the Cartan relations for the
calculus defined by $\CM_0=km_0=\sum\limits_{v\in cent(t_0)}
k(v.t_0)$ and $C_0^X.$ We label the basis of $\CM$ using elements
of $C_0^X$ as
\[f_{z_{ij}}:=m_0\tra \overline{z_{ij}}=\sum\limits_{v\in cent(t_0)}
\und{\bar {s_j}}^{-1}v\und {\bar {s_i}}.s_j^{-1}\] and then denote
by $(e_{z_{ij}})$ the dual basis of $(f_{z_{ij}})$.
\begin{lemma}\label{adjointofaction}
The  action $*$ on the basis $(e_{z_{ij}})$ is
\[e_{z_{lm}}*(e.s_j)=e_{s_j^{-1}z_{lm}s_j},\quad
 e_{z_{lm}}*(u.e)=e_{u^{-1}{z_{lm}}u},\quad u,s_j\in G\]
\end{lemma}
for all $0\leq j,m,l\leq N$, i.e. $X$ acts by the right adjoint action on the indexes.\\
\proof From the definition of $\tra$ we have for all $0\leq
p,q,j\leq N$
\begin{eqnarray}\label{actionf_z_pq}
f_{z_{pq}}\tra s_j^{-1}=f_{z_{pq}}\tra
e.s_j^{-1}&=&\sum\limits_{v\in cent(t_0)} (\und {\bar {s_q}}^{-1}v
\und {\bar {s_p}}.s_q^{-1})
\tra s_j^{-1}\nonumber\\
&=&\sum\limits_{v\in cent(t_0)}{\und {s_j}} \und {\bar
{s_q}^{-1}}v \und {\bar {s_p}} \und {s_j^{-1}}.
s_j{s_q}^{-1}s_j^{-1}
\end{eqnarray}
On the other hand, $f_{z_{pq}}\tra s_j^{-1}$ is homogeneous and
should be
 linear combination of $f_{z_{ij}}, 0\leq i,j\leq N.$ But the latter have
different degrees then we deduce that $f_{z_{pq}}\tra s_j^{-1}$ is
linear combination of only one of them, the one  whose degree is
$||f_{z_{pq}}\tra s_j^{-1}|| = s_jz_{pq}s_j^{-1}$, explicitly,
$f_{z_{pq}}\tra s_j^{-1}=cf_{s_jz_{pq}s_j^{-1}}$ where $c$ is a
constant. In fact this constant is 1 since in the expansion of
$f_{z_{pq}}\tra s_j^{-1}$ in equation (\ref{actionf_z_pq}) the
nonzero coefficients of $us \in X$ equal 1. Therefore
$f_{z_{pq}}\tra s_j^{-1}=f_{s_jz_{pq}s_j^{-1}}$ and
\begin{eqnarray*}
e_{z_{lm}}*(e.s_j)&=&
\sum\limits_{p,q}<e_{z_{lm}},f_{z_{pq}}\tra s_j^{-1}>e_{z_{pq}}\\
&=&\sum\limits_{p,q}<e_{z_{lm}},f_{s_jz_{pq}s_j^{-1}}>e_{z_{pq}}\\
&=&e_{s_j^{-1}z_{lm}s_j}
\end{eqnarray*}
as stated. One  follows the same reasoning  to prove  the second
assertion of the lemma. \eproof

We can now explicitly give the differential calculus of dimension
$|C^G|^2$ defined by $(m_0,C^G)$ as above for each conjugacy class
$C^G$ of $G$. We use Theorem~\ref{cartanrelations} and
 Proposition~\ref{braidingonVXV}.
\begin{prop}\label{canonicalcalculusonD(G)}
The Cartan calculus and braiding for the canonical differential
calculus on $D^*(G)$ defined by $(m_0,C_0^G)$ are given by:\\
(i)\quad Commutations relations
\[e_{z_{ij}} f=R_i(f)e_{z_{ij}},
\quad e_{z_{ij}} u=(\und {s_j}^{-1}u\und {s_j})e_{u^{-1}z_{ij}u}\]
(ii)\quad Differentials
\[df=\sum\limits_i\partial_i(f)e_{z_{ii}},
\quad du=\sum\limits_i(\und {s_i}^{-1}u \und
{s_i})e_{u^{-1}z_{ii}u}-\sum\limits_iue_{z_{ii}}\] (iii)\quad The
element
 \[\theta=\sum\limits_ie_{z_{ii}}\]
(iv)\quad The braiding
\[\Psi(e_{z_{ij}}\otimes e_{z_{lm}})=e_{s_j^{-1}z_{lm}s_j}
\otimes e_{\und{s_l}^{-1}\und{s_m}z_{ij}\und{s_m}^{-1}\und{s_l}}\]
for $u\in G$, $f \in k(M)=k(G)$, where $R_i(f)(g)=f(gs_i^{-1})$
for all $g\in G$, and $\partial_i=R_i-\id $.
\end{prop}
\proof Since $$||f_{z_{ij}}||=z_{ij}=\und {s_i}\und
{s_j}^{-1}.s_j,$$ we have \[\<f_{z_{ij}}\>=s_j^{-1}\ra
(s_is_j^{-1})^{-1}=e.s_j^{-1},\quad |f_{z_{ij}}|^{-1}=s_j^{-1}\la
(\und {s_i}\und {s_j}^{-1})^{-1}=\und{s_i}^{-1}\und{s_j}.e\] then
\begin{eqnarray*}
<\delta_{\<f_{z_{ij}}\>},f_{z_{ij}}>&=&\sum\limits_{v\in
cent(t_0)} <\delta_{s_j^{-1}},\und {\bar {s_j}}^{-1}v
\und {\bar {s_i}}.s_j^{-1}>\\
&=&\sum\limits_{v\in cent(t_0)}\delta_{v\und {\bar {s_i}}, \und
{\bar {s_j}}}\ =\ \delta_{\bar {s_i},\bar {s_j}}=\delta_{i,j}
\end{eqnarray*}
where we use the fact that for  $ v\in cent(t_0),$ $\bar
{s_j}=v\bar {s_i} \Longrightarrow \bar {s_j}^{-1}t_0^{-1}\bar
{s_j}=\bar {s_i}^{-1}t_0^{-1}\bar {s_i}=s_i$ and  by definition of
$\bar {s_i}$ we deduce $\bar {s_j}=\bar {s_i}$. We then rewrite
the results in Theorem~\ref{cartanrelations}, Corollary
\ref{theta} and Proposition~\ref{braidingonVXV} using the previous
Lemma~\ref{adjointofaction} the $G-M$ bigrading and pairing above
to obtain the results as stated. \eproof

One may verify that the restriction to $k(G)$ of the differential
calculus $(m_0,C_{t_0^{-1}})$ on $D^*(G)$ is exactly the
differential calculus defined on $k(G)$ by $C_{t_0}$ in
Proposition~\ref{groupcase}(i) after suitable matching of the
conventions. These results from our theory for bicrossproducts are
in agreement with calculi on $D^*(G)$ that can be constructed by
entirely different methods \cite{MAj1} via its coquasitriangular
structure.

\subsection{The case $X=G\lcross M$ with $G$ Abelian}
It is known \cite{MAjFounda} that if $G$ is Abelian then $kG\cong
k(\hat G)$ and equivalently $k(G)\cong k\hat{G},$ where $\hat G$
is the group of characters of $G$. Then
$$A=k(M)\lcross kG=k(M)\lcross k(\hat G)\cong k(M\rcross\hat G).$$
The product in $M\rcross\hat G$ is
$$(t.\psi)(s.\phi)=(ts. (\psi\ra s)\phi)$$
where we denote the  element $(t,\psi)$ by $t.\psi$, using
factorization notation. The action of $M$ on $\hat G$ is
$$(\psi\ra s)(u)=\psi(s\la u), \quad \forall s \in M, \forall u \in
G.$$ Explicitly an element $f\in k(G)$ is viewed as
\[\tilde f=\sum\limits_{\phi \in \hat G,u\in G}\frac{1}{|G|}\phi (u^{-1})f(u)\phi \in k\hat G\]
while $\phi \in \hat G$ is viewed as
\[\tilde \phi =\sum\limits_{u\in G}\phi (u)\delta_u \in k(G)\]
This induces Hopf algebras (Fourier) isomorphisms
\[\CF:k(M)\lrbicross kG)\longrightarrow k(M\rcross \hat G),\quad
\CF^*:k.M\rcross \hat G\longrightarrow  kM\rcross k(G)\]

defined by

\[ \CF(\delta_t \tens v))=\sum\limits_{\chi \in \hat G}\chi
(v)\delta_{t.\chi},\quad
\CF^{-1}(\delta_{t.\chi}))=\sum\limits_{u\in
  G}\frac{1}{|G|}\chi(u^{-1})\delta_t\tens u\]
and
 \[\CF^*(s.\chi)=\sum\limits_{u\in G}\chi(u)s\tens \delta _u, \quad
{\CF^*}^{-1} (s\tens \delta_u)=\sum\limits_{\chi \in \hat
G}\frac{1}{|G|}\chi (u^{-1} )s.\chi\] $\forall t,s \in M, \forall
u\in G,\forall \chi \in \hat G.$

Since $A$ is isomorphic to the algebra of functions on a group, it
follows that the irreducible bicovariant differential calculi on
$A$ from the general theory above must correspond to nontrivial
conjugacy classes of $M\rcross\hat G$. We now exhibit this
correspondence as follows:

For the first direction, let $\hat C_0$ be a nontrivial  conjugacy
class of $t_0.\psi _0$ in $M\rcross\hat G.$ This class defines an
irreducible bicovariant differential calculus on $A=
k(M\rcross\hat G)$ whose quantum tangent space is
$$L=k.\{a-e, a\in \hat C_0\}\subset
\ker\varepsilon\subset A^*=kM\rcross k(G).$$ From this we
determine $\CM\subset X$ as
$$k\bar 1\oplus \CM=\Pi ^{-1}(L)\quad \hbox{(by Theorem \ref{tangentspacestheo})}$$
Then we take as conjugacy class $C^X$ in $Z$ that determined by
the conjugacy class $C_0^M$ of $t_0^{-1}$ in $M$, namely $C_0^X$,
and we define $\CM_0$ using Proposition~\ref{M_0fromM} as
$$\CM_0=\CM\tra \delta_{t_0^{-1}}.$$ One then verifies easily that
$\CM\tra \delta_{t_0^{-1}}$ is
nonzero as required in Lemma \ref{M_0fromM} and  $C^X$ defined
above does not depend on the chosen element in $\hat C_0$.

For the second direction, we suppose that we are given a nonzero
irreducible bicovariant differential calculus on $A$ defined (say)
by an irreducible subrepresentation $\CM\subset kX$ under the
action of $D(X)$. We need to construct a conjugacy class $\hat C
\subset M\rcross\hat G$   such that the differential calculus
defined on $A$ by $\hat C$ coincides with that defined by $\CM$,
i.e.,
\begin{eqnarray}
k.(\hat C-e):=k.\{a-e, a\in \hat C\}=\Pi(\CM)
\end{eqnarray}
as quantum tangent spaces in $H=A^*$. First of all, we note that
$H$ is the group algebra $k. M\rcross\hat G$ so that
$\ker\varepsilon _H$ is generated as vector space by the set
$$B_{\varepsilon _H}=\{t.\psi-e, \quad t\in M, \psi \in \hat G\}.$$ Since
$\Pi(\CM) \subset \ker\varepsilon _H$, for all $m\in \CM$,
$\Pi(m)$ is linear combination of elements of $B_{\varepsilon
_H}$. In general, not all of such elements are necessary to span
$\Pi(\CM)$, so let us denote by
$$B_{\CM}=\{t_i.\psi_j-e, \quad t_i\in M, \psi_j \in \hat G,
\quad (i,j)\in I\times J\}$$ a minimal set of elements of
$B_{\varepsilon _H}$ such that
$$\Pi(\CM)\subset k.B_{\CM}.$$
A long but not difficult computation using Fourier isomorphisms
above shows that $k.B_{\CM}=\Pi(\CM)$ and any conjugacy class
$\hat C$ in $M\rcross\hat G$,  of an element $t_1.\psi_1$ such
that $t_1.\psi_1-e\in \Pi(\CM)$  obeys
$$k.(\hat C-e)=\Pi(\CM)$$
as expected.

\section{Canonical calculi on crossproducts
 $k(M)\lcross kG$}\label{crossproductcase}

We now consider the complementary special case where $X=G\rcross
M$ and $A=k(M)\lcross kG,$ a cross product Hopf algebra. We show
that conjugacy classes in $M$ which are invariant under the right
action of $G$ define canonical bicovariant differential calculi on
both $k(M)$ and on $A$ such that the calculus on $A$ is an
extension of the one on $k(M)$. This gives a natural way to define
bicovariant differential calculi on the double $D(G)$ of any
finite group $G$.

\begin{prop}\label{extendedcalonA}
Let $X=G\rcross M$ be a semi-direct factorization. For any
$G$-invariant conjugacy class $C$ of $ M$, when it exists, we set
\begin{eqnarray}
\CM=\bigoplus_{a\in C}k(\sum\limits_{v\in G}v.(a^{-1}\ra
v^{-1}))\subset kX
\end{eqnarray}
Then\\
(i)\quad The vector space $\CM$ is isomorphic to an irreducible
quantum tangent space in $kM\rcrossco k(G)$, precisely it is of
the form
$\CM=\bigoplus_{z\in C} \CM_0 \tra \bar z$
as above.
\\
(ii)\quad The differential calculus defined  on  $kM\rcross k(G)$
by $(\CM, C)$ restricts to the  calculus defined on $k(G)$ by
$C^{-1}$ as in
 Proposition~\ref{groupcase}.

We call the calculus $(\CM, C)$, the canonical differential
calculus defined on $kM\rcross k(G)$ by the conjugacy class $C$.
\end{prop}

\proof We now have $Z=M$ since $v^{-1}t^{-1}v= t^{-1}\ra v$ for
all $t\in M, v\in G.$ and the action in (\ref{actionD(X)onkX})
becomes
\[vt\tra us=vu.\tilde s t {\tilde s}^{-1},\quad \tilde s =s^{-1}\ra (vu)^{-1}.\]
Let $C$ be a nontrivial  $G$-invariant conjugacy class in $M$ and
let $t_0\in C$ and $X_{t_0}$ the centralizer of $t_0$. then
\[\eta ^{-1}(t_0)=\{vt\in X,t^{-1}\ra v=t_0\}=\{v(t_0^{-1}\ra
v^{-1}),\;v\in G\}\] Let $us\in X_{t_0}$ we have
\[ust_0=t_0us\Longrightarrow t_0^{-1}\ra u^{-1} =ut_0^{-1}u^{-1}
=(s^{-1}\ra u^{-1})t_0^{-1}(s\ra u^{-1})\]
so that the action of  $us\in X_{t_0}$ on $v(t_0^{-1}\ra v^{-1})
\in \eta ^{-1}(t_0)$ is
\begin{eqnarray*}
v(t_0^{-1}\ra v^{-1}) \tra us&=&vu. ((s^{-1}\ra u^{-1})\ra v^{-1})
(t_0^{-1}\ra v^{-1})((s\ra u^{-1})\ra v^{-1})\\
&=&vu.((s^{-1}\ra u^{-1})t_0^{-1}(s\ra u^{-1})\ra v^{-1})\\
&=& vu.(t_0^{-1}\ra (vu)^{-1})\in \eta ^{-1}(t_0)
\end{eqnarray*}
This proves that the one-dimensional vector space
\[\CM_{t_0} =k.\sum\limits_{v\in G}v.(t_0^{-1}\ra v^{-1})\]
is an irreducible $X_{t_0}$-module.

For any $a\in C\subset M$, we fix $\bar a \in M$ such that
$a={\bar
  a}^{-1}t_0\bar a$ and set
\[\CM=\bigoplus_{a\in C} \CM_{t_0} \tra \bar a.\]
It is clear that the G-invariance of $C$ implies that $C$ is also
a conjugacy class in $X$ hence by
Proposition~\ref{centralizermodule} and
Theorem~\ref{tangentspacestheo}, $\CM$ is isomorphic to a quantum
tangent space in $kM\rcrossco k(G).$

Furthermore we easily obtain
\begin{eqnarray*}
\CM_{t_0} \tra  \bar a= k\sum\limits_{v\in G}v.({\bar
  a}^{-1}t_0^{-1}\bar a
\ra v^{-1})=k.\sum\limits_{v\in G}v(a^{-1}\ra v^{-1})
\end{eqnarray*}
then
\[\CM=\bigoplus_{a\in C}k(\sum\limits_{v\in G}v.(a^{-1}\ra v^{-1})).\]

For the Cartan calculus of the differential calculus associated to
$\CM$, we choose the canonical basis $(f_a)_{a\in C}$ defined by
\[f_a=:\sum\limits_{v\in G}v.(a^{-1}\ra v^{-1}) \]
then it is clear that
\[\<f_a\>=a^{-1},\quad |f_a|=e,\quad <\delta_{\<f_a\>},f_a>=1,\quad
\forall a\in C.\] On the other hand, for $u\in G$ and $c\in C,$
\[ f_c\tra u^{-1}=\sum\limits_{v\in G}vu^{-1}.(c^{-1}\ra v^{-1})=
\sum\limits_{w\in G}w.((c^{-1}\ra u^{-1})\ra w^{-1})=f_{c^{-1}\ra
  u^{-1}}=f_{ucu^{-1}}\]
so that
\[e_a*u=\sum\limits_{c\in C}<e_a, f_c\tra u^{-1}>e_c=e_{u^{-1}au}.\]

Then the Cartan calculus from Theorem~\ref{cartanrelations} reads
\[e_a\delta_s=\delta_{sa}e_a,\quad e_au=ue_{u^{-1}au}\]
\[d\delta_s=\sum\limits_{a\in C} (\delta_{sa}-\delta_s)e_a,\quad
du=\sum\limits_{a\in C}u(e_{u^{-1}au}-e_a)\] It is now clear that
this canonical differential calculus on $k(M)\lcross kG$ is an
extension of the differential calculus defined on $k(M)$ as in
Proposition~\ref{groupcase} by the opposite conjugacy class
$C^{-1}$ of $C$. \eproof

\section{Exterior algebra and cohomology computations}

We have already seen that $\theta$ in Corollary~\ref{thetainv} is
a nontrivial element of the noncommutative de Rham cohomology for
any bicrossproduct. In this section we will glean more insight
into the cohomology through a close look at particular
bicrossproducts. From a physical point of view this is the
beginning of `electromagnetism' on such spaces. From a
noncommutative geometers point of view it is the `differential
topology' of the algebra equipped with the differential structure.
Note that very little is known in general about the full
noncommutative de Rham cohomology even for finite groups, but
insight has been gained through examples such as in
\cite{MAjRai,ngaMAj}. We are extending this process here.

In particular, just as for a Lie algebra there is a unique
differential structure giving a connected and simply connected Lie
group, so we might hope for a `natural' if not unique choice of
calculus such that at least $H^0=k.1$, which is a connectedness
condition (so that a constant function is a multiple of the
identity) and with small $H^1$. By looking at several examples and
using our explicit Cartan relations for bicrossproducts, we find
that a phenomenon of this type does appear to hold.  In
particular, as a main result of the paper from a practical point
of view, we find a unique such calculus on the quantum double
$D(S_3)$ viewed as a bicrossproduct, i.e. a natural choice for its
differential geometry. We also cover the codouble $D^*(S_3)$ as
another bicrossproduct.

In each case studied here,  we describe the factorizing groups,
the set $Z$ and hence the classification of calculi. We then
compute the first order calculi in each case using the theory
above, and the braiding on basic forms $\{e_a\}$ dual to the basis
$\{f_a\}$ stated in each case of the quantum tangent space yielded
by the classification. In each case,
\[\Omega^1(A)=A\tens \Lambda^1,\quad \Lambda^1:=\<e_a\>_k\]
where $A$ is the bicrossproduct Hopf algebra and $<\ >_k$ denotes
the $k$-span. In describing the exterior derivative we use the
translation and `finite difference' operators
\[ R_s(f)(r)=f(rs);\quad (\del_s f)(r)=f(rs)-f(r),\quad
\forall r\in M,\quad f\in k(M)\] for the relevant group $M$ and
relevant $s\in M$, as already used elsewhere.

From the braiding we then compute the higher order differential
calculus using the braided factorial matrices $A_n$ given by
\[A_n= ({\rm id}\otimes A_{n-1})[n,-\Psi],\; [n,-\Psi]
={\id}-\Psi_{12}+
\Psi_{12}\Psi_{23}+...+(-1)^{(n-1)}\Psi_{12}...\Psi_{n-1,n}\]
where $\Psi_{i,i+1}$ denotes $\Psi$ acting in the $i,i+1$
positions in $\Lambda^1{}^{\tens n}$. The space $\Lambda^n$ of
invariant $n$-forms is then the quotient of $(\Lambda^1)^{\otimes
n}$ by $\ker A_n$. This is the computationally efficient braided
groups approach used in \cite{MAjRai,MAjYangmills,ngaMAj} and
equivalent to the original Woronowicz description of the
antisymmetrizers in \cite{W}. These braided integer matrices have
also been adopted by other authors, such as \cite{Ros}.

\subsection{Calculi and cohomology on $k(\Z_2)\rcrossco k\Z_3$}

This baby example $k(\Z_2)\rcrossco k\Z_3$ is actually a
semidirect coproduct isomorphic to $k(S_3)$ and among other things
demonstrates the Fourier theory in Section~4.2. From the theory of
calculi on finite groups, we know that there are two irreducible
calculi of dimensions 2,3 respectively, according to the
nontrivial conjugacy classes of $S_3$. We illustrate how this
known result comes about in our bicrossproduct theory.

Here, $X=S_3$ factorizes into $M=\Z_2=\{e,s\}$ and
$G=\Z_3=\{e,u,u^2\}$, where $s=(12),u=(123)$. The right action of
$G$ on $M$ is trivial and the left action of $M$ on $G$ is defined
by $s\la =(u,u^2)$ (the permutation). The set $Z$ of elements
$||x||$ is $Z=\{e,s,us,u^2s\}$ which splits into two conjugacy
classes $C^X=\{e\}$ and $C^X=\{s,us,u^2s\}$. This leads to the
following irreducible bicovariant calculi.

\vspace{0.5cm}

$(i)\quad C^X=\{e\},$ $\CM=<f_1,f_2>_k$, where $q=e^{2\pi i\over
3}$ and
\[f_1=e+q^2 u+q u^2,\quad f_2=e+q u+q^2 u^2\]
\[e_af=f e_a,\quad  \forall f\in k(\Z_2);\quad  a=1,2\]
\[e_1u=q^2 ue_1,\quad e_1u^2=q u^2e_1,\quad e_2u=q ue_2,
\quad e_2u^2=q^2u^2e_2\]
\[\theta=e_1+e_2\]
\[ df=0,\ \forall f\in k(\Z_2),\quad du=u(q^2-1)e_1+u(q-1)e_2,
\quad du^2=u^2(q-1)e_1+u^2(q^2-1)e_2\]
\[\Psi(e_a\otimes e_b)=e_b\otimes e_a,\ a,b=1,2.\]\\
The exterior algebra has the usual relations and
dimensions\[e_a^2=0,\quad e_1\wedge e_2=-e_2\wedge e_1,\quad
\dim(\Omega)= 1:2:1.\] The cohomology can be identified with \[
H^0=k(\Z_2),\quad H^1=k(\Z_2)e_1\oplus k(\Z_2) e_2,\quad
H^2=k(\Z_2)e_1\wedge e_2\] with dimensions 2:4:2.

\vspace{0.5cm}
 $(ii)\quad C^X=\{s,us,u^2s\}$, $\CM=<f_1,f_2,f_3>_k$, where
\[f_1=s,\quad f_2=us,\quad f_3=u^2s\]
\[e_a f=R_s(f) e_a,\quad \forall f\in k(\Z_2),\quad  a=1,2,3\]
\[e_1u=u^2 e_2,\ e_1u^2= ue_3,\ e_2u= u^2e_3,\ e_2u^2=ue_1,
\ e_3 u=u^2e_1,\ e_3u^2=ue_2\]
\[\theta=e_1,\ d f=\del_s(f)e_1,\quad du=u^2e_2-ue_1,
\quad du^2=ue_3-u^2e_1\]
\begin{eqnarray*}
\Psi(e_1\otimes e_1)&=&e_1\otimes e_1,\quad \Psi(e_2\otimes
e_1)=e_1\otimes e_2,
\quad \Psi(e_3\otimes e_1)=e_1\otimes e_3\\
\Psi(e_1\otimes e_2)&=&e_3\otimes e_3,\quad \Psi(e_2\otimes
e_2)=e_3\otimes e_1,
\quad \Psi(e_3\otimes e_2)=e_3\otimes e_2\\
\Psi(e_1\otimes e_3)&=&e_2\otimes e_2,\quad \Psi(e_2\otimes
e_3)=e_2\otimes e_3,
\quad \Psi(e_3\otimes e_3)=e_2\otimes e_1\\
\end{eqnarray*}
The exterior algebra is quadratic with relations \[ e_1\wedge
e_1=0,\quad e_2\wedge e_3=0,\; e_3\wedge e_2=0,\quad e_1\wedge
e_2+e_2\wedge e_1+e_3^2=0,\quad e_1\wedge e_3+e_3\wedge
e_1+e_2^2=0\] and has dimensions and cohomology:
\[\dim(\Omega)=1:3:4:3:1\]
\[ H^0=k.1,\quad H^1=k.\theta,\quad H^2=0,\quad
H^3=ke_3^3,\quad H^4=k.e_1\wedge e_2\wedge e_2 \wedge e_2.\]

This is isomorphic to the cohomology and calculus on $k(S_3)$
studied in \cite{MAjRai}. We see that this is the unique choice
$(ii)$ with $H^0=k.1$. We also see that both calculi exhibit
Poincar\'e duality.

\subsection{Calculi and cohomology on $k(S_3)\lrbicross k\Z_6$}

The second example  $k(S_3)\lrbicross k\Z_6$ is a nontrivial
bicrossproduct\cite{BGM} but is (nontrivially) isomorphic to a
version of the dual of a quantum double $D^*(S_3)=k(S_3)\rcrossco
kS_3$. Among other things, it demonstrates our results for the
codouble in Section~4.1.

Here $X=S^3\times S^3$ factorizes differently into groups
\[ G=\Z_6=\{u^0,u,u^2,u^3,u^4,u^5\},\quad
M=S_3=\{e,s,t,t^2,st,st^2\}\] where $s=(12), t=(123)$  and $u$ is
the generator of $\Z_6$. The right action of $u$ on $M$ is the
permutation \[\ra u=(st,st^2)(t,t^2)\] while the left action of
$M$ on $G$ is given completely in terms of permutations  by
\begin{eqnarray*}
e\la&=&\id,\quad s\la = (u,u^5)(u^2,u^4),\quad
t\la=(u^5,u^3,u),\quad
 t^2\la=(u,u^3, u^5)\\
st\la &=&(u^2,u^4)(u^3, u^5),\quad st^2\la=(u^2,u^4)(u,u^3)
\end{eqnarray*}
For $X=\Z_6.S_3$ the set of the values of $||.||$ is
\[Z=\{e,s,t,t^2,st,st^2,u^2s,u^4s,u^2t,u^4t^2,u^2st,u^4st,
u^2st^2,u^4st^2\}\] which splits into three conjugacy classes
\begin{eqnarray*}
C^X=\{0\},\quad C^X=\{t,t^2,u^2t,u^4t^2\},\quad
C^X=\{s,st,st^2,u^2s,u^4s,u^2st,u^4st,u^2st^2,u^4st^2\}.
\end{eqnarray*}
If we choose the respective basis points to be $z_0=e,t,s$ then we
have
\[\CN^{-1}(e)=G,\quad \CN^{-1}(t)=\{t^2,u^2t^2,u^4t^2\},\quad
 \CN^{-1}(s)=\{s,u^3s\}\]
The centralizers of $e,t,s$ in $X$ are respectively
\[X_e=X,\quad G_t=\{e,t,t^2,u^2,u^2t,u^2t^2,u^4,u^4t,u^4t^2\},\quad
G_s=\{e,s,u^3,u^3s\}\] Applying the general theory of Sections 4
and 5 to  these data leads to the irreducible bicovariant
differential calculi on $A=k(S_3)\lrbicross k\Z_6$ as follows:
\vspace{0.4cm}

 $(i)\quad C^X=\{e\}$, $\CM=<f_1>_k$, where
$f_1=u^0-u+u^2-u^3+u^4-u^5$.
\[e_1f=fe_1,\quad \forall f\in k(S_3),\quad e_1u^i=(-1)^iu^ie_1,
\forall u^i\in G\]
\[\theta=e_1,\quad df=0,\forall f\in k(S_3),
\quad  du^i=(-1+(-1)^i)u^ie_1\]
\[ \Psi(e_1\otimes e_1)=e_1\otimes e_1\]

\vspace{.4cm}

$(ii)\quad C^X=\{e\}$, $\CM=<f_1,f_2,f_3,f_4>_k$,where setting
$q=e^{-2\pi i\over 6}$,
\[f_1=u^0+q u +q^2 u^2-u^3-q u^4-q^2 u^5,\quad
f_2=u^0+q^2 u -q u^2+u^3+q^2 u^4-q u^5\]
\[f_3=u^0-q u +q^2 u^2+u^3-q u^4+q^2 u^5,\quad
f_4=u^0-q^2 u -q u^2-u^3+q^2 u^4+q u^5\]
\[e_a f=f e_a,\quad \forall f\in k(S_3),\quad a=1,2,3,4,\quad e_1u^j=q^ju^je_1,
\;\;e_2u^j=q^{2j}u^je_2,\]
\[
e_3u^j=q^{4j}u^je_3,\quad e_4u^j=q^{5j}u^je_4,\quad \theta=e_1+e_2+e_3+e_4\]
\[d f=0,\forall f\in k(S_3),\quad du^j=u^j((q^j-1)e_1
+(q^{2j}-1)e_2+(q^{4j}-1)e_3+(q^{5j}-1)e_4)\]
\[ \Psi(e_a\otimes e_b)=e_b\otimes e_a,\quad a,b=1,2,3,4\]

\vspace{.4cm}

$(iii)-(v)\quad C^X=\{t,t^2,u^2t,u^4t^2\}$, $
\CM=<f_1,f_2,f_3,f_4>_k$, where $q=1,e^{2\pi\imath\over
3},e^{-2\pi\imath\over 3}$ for the three cases and
\[f_1=t^2+q u^2t^2+q^2 u^4t^2,\quad f_2=t+q^2 u^2t+q u^4t,\]
\[f_3=ut^2+q u^3t^2+q^2 u^5t^2,\quad f_4=u^5t+q^2 ut+q u^3t\]
\[e_1f=R_t(f)e_1,\quad e_2f=R_{t^2}(f)e_2,
\quad e_3f=R_{t^2}(f)e_3,\quad e_4f=R_t(f)e_4,\quad \forall f \in
k(S_3)\]
\begin{eqnarray*}
e_1u^{2j}&=&q^{-2j}u^{2j}e_1,\quad e_1u^{2j+1}=q^{-2j}u^{2j+3}e_3,\\
e_2u^{2j}&=&q^{-j}u^{2j}e_2,\quad e_2u^{2j+1}=q^{2-j}u^{2j+5}e_4\\
e_3u^{2j}&=&q^{-2j}u^{2j}e_3,\quad e_3u^{2j+1}=q^{1-2j}u^{2j+5}e_1,\\
 e_4u^{2j}&=&q^{-j}u^{2j}e_4,\quad e_4u^{2j+1}=q^{-j}u^{2j+3}e_2
\end{eqnarray*}
\[\theta=e_1+e_2,\quad \extd f=\del_t(f)e_1+\del_{t^2}(f)e_2\]
\[ du^{2j}=(q^{-2j}-1)u^{2j}e_1+(q^{-j}-1)u^{2j}e_2,\]
\[du^{2j+1}=q^{-2j}u^{2j+3}e_3+q^{2-j}u^{2j+5}e_4-u^{2j+1}(e_1+e_2)\]
\begin{eqnarray*}
\Psi(e_a\otimes e_1)&=&e_1\otimes e_a,\quad \Psi(e_a\otimes e_2)
=e_2\otimes e_a,\quad a=1,2,3,4\\
\Psi(e_1\otimes e_3)&=&q e_3\otimes e_1,\quad \Psi(e_1\otimes
e_4)=q^2 e_4\otimes e_1,
\quad \Psi(e_2\otimes e_3)=q^2 e_3\otimes e_2,\\
  \Psi(e_2\otimes e_4)&=&q e_4\otimes e_2,
  \quad  \Psi(e_3\otimes e_3)= e_3\otimes e_3,\quad
 \Psi(e_3\otimes e_4)= e_4\otimes e_3,\\
\Psi(e_4\otimes e_3)&=& e_3\otimes e_4,\quad \Psi(e_4\otimes e_4)=
e_4\otimes e_4
\end{eqnarray*}

The resulting exterior algebra and cohomology depend on the
braiding. In case $(iii)$ we have:
\[ \dim(\Omega)=1:4:6:4:1,\quad H^0=<\delta_{s^i} u^{2j}|\ i=0,1;\
j=0,1,2>_k\] \[H^1=<\delta_{s^i} u^{2j}e_a|\ i=0,1;\ j=0,1,2;\quad
a=1,2,3,4>_k\] where the cohomology is 6-dimensional in degree 0
and 24 dimensional in degree 1. The relations in the exterior
algebra are that the forms $\{e_a\}$ anticommute as usual. In case
$(iv)$ we have the 6 relations
\[ e_a^2=0,\quad e_1\wedge e_2+e_2\wedge e_1=0,
\quad e_3\wedge e_4+e_4 \wedge e_3=0\] and
\[ \dim(\Omega)=1:4:10:53:\cdots,
\quad H^0=<\delta_{s^i} \phi_j u^{2j}|\ i=0,1;\ j=0,1,2>_k\] \[
H^1=<\delta_{s^i} \phi_j u^{2j}e_a|\ i=0,1;\ j=0,1,2;\quad
a=1,2>_k\] where
\[ \phi_i=\sum_{j=0}^{j=2}\delta_{t^j}q^{ij}.\]
 Here the dimensions of the cohomology are 6 in degree 0
and 12 in degree 1. The case $(v)$ is identical with $q$ replaced
by $q^{-1}$.

\vspace{.5cm} $(vi)-(vii)$ \quad
$C^X=\{s,st,st^2,u^2s,u^4s,u^2st,u^4st,u^2st^2,u^4st^2\}$,
\[ \CM=<f_{ai}|\ a=1,2,3;\ i\in\Z_3>_k\]
where $q=\pm1$ for the two cases and:
\[f_{10}=s+qu^3s,\  f_{11}=u^2s+qu^5s, \
 f_{12} =us+qu^4s,\  f_{20}=st+qu^5st^2,\ f_{21}=
 u^2st+qust^2\]
\[ f_{22}=u^4st+qu^3st^2,\quad f_{30}=st^2+qust,\quad
f_{31}=u^3st+qu^2st^2,\quad f_{32}=u^5st+qu^4st^2.\] For brevity,
we give the details only for the $q=1$ case (the other is
similar). Then
\[ e_{1i}f=R_s(f)e_{1i},\quad e_{2i}f=R_{st}(f)e_{2i},\quad e_{3i}f
=R_{st^2}(f)e_{3i},\quad \forall f\in k(S_3)\]
\[e_{1i}u^j=u^{-j}e_{1,i-j},\;\;  e_{2i}u^{2k}=u^{4k}e_{2,i+k},\;\;
e_{2i}u^{2k+1}=u^{4k+1}e_{3,i+k}\]
\[e_{3i}u^{2k}=u^{4k}e_{3,i+k},\quad
e_{3i}u^{2k+1}=u^{4k+3}e_{2,i+k+1}\]
\[\theta=e_{10}+e_{20}+e_{30},\quad d f=\del_s(f)e_{10}
+\del_{st}(f)e_{20}+\del_{st^2}(f)e_{30}\]
\[ du^{2k}=u^{-2k}e_{1,k}+u^{4k}e_{2,k}+u^{4k}e_{3,k}-u^{2k}\theta \]
\[du^{2k+1}=u^{-(2k+1)}e_{1,-(2k+1)}+u^{4k+3}e_{2,k+1}
+u^{4k+1}e_{3,k}-u^{2k+1}\theta\]
\[ \Psi(e_{1i}\tens e_{aj})=e_{(23)a,-j}\tens
e_{1,i-j}\]
\[\Psi(e_{2i}\tens e_{aj})=e_{(13)a,-j}\tens
e_{2,i-j},\quad \Psi(e_{3i}\tens e_{aj})=e_{(12)a,-j}\tens
e_{3,i-j}.\]

The resulting exterior algebra has relations
\[ e_{aj}\wedge e_{a,-j}=0,\quad
e_{ai}^2+\{e_{a,i-1}, e_{a,i+1}\}=0\]
\[ e_{1i}\wedge e_{2j}+e_{2,j-i}\wedge e_{3,-i}+e_{3,-j}
\wedge e_{1,i-j}=0\]
\[
e_{2j}\wedge e_{1i}+e_{3,-i}\wedge e_{2,j-i}+e_{1,i-j}\wedge
e_{3,-j}=0\] for $a=1,2,3$ and $i,j\in\Z_3$. The dimensions of the
exterior algebra and cohomology in low degree are
\[ \dim(\Omega)=1:9:48:198:\cdots,\quad H^0=k.1,\quad
H^1=k.\theta.\]

From these explicit computations we conclude in particular:
\begin{prop} Only
the 9-dimensional calculi $(vi)$--$(vii)$ have $H^0=k.1$
\end{prop}

The natural one here is $(vi)$ where $q=1$ with the other as a
signed variant. We also have Poincar\'e duality at least for all
cases where the exterior algebra was small enough to be fully
computed. According to \cite{BM} this bicrossproduct is a
coquasitriangular Hopf algebra, isomorphic  to the quantum
codouble $D^*(k(\Z_2)\rcrossco k\Z_3)$ of our first example in
Section~6.1, hence also to the quantum codouble
$D^*(S_3)=k(S_3)\rcrossco kS_3$ of the type covered in
Section~4.1. The canonical calculi given in
Theorem~\ref{canonicalcalculusonD(G)} correspond to $(iii)$ and
$(vi)$, with $(vi)$ indeed the canonical extension of the natural
(3-dimensional) calculus on $S_3$ (as in Section~6.1). The other
cases fit in their number and dimensions with a completely
different classification theorem for factorizable
coquasitriangular Hopf algebras\cite{MAj1} which implies that
calculi can be classified by representations of the quantum double
$D(S_3)$, with dimension the square of that of the representation.
These are labelled by conjugacy classes in $S_3$ and
representations of the centralizer, giving calculi of dimensions
1,4 and 4,4,4 and 9,9 for the three classes. We see that we obtain
isomorphic results from our bicrossproduct classification (the
isomorphism is nontrivial, however).

\subsection{Calculus and cohomology on $D(S_3)$}

Finally we consider the dual example to the preceding one, with
``coordinate ring'' the cross product
 $k(S_3)\lcross kS_3 =D(S_3)$ in the same conventions as for $D(X)$
in Section~\ref{moduleofD(H)}. This corresponds to the semidirect
factorization in Section~\ref{crossproductcase}, namely
$X=S_3\rcross S_3$ with action by conjugation.

Here $G=M=S_3$, thus we use the same notation as
 in Section~\ref{Gactstrivaily} namely elements of $S_3.e$ are underlined,
those of $e.S_3$ are not, so that a general element of $X=S_3.S_3$
is of the form  $\und{v}.t$. We set again
$S_3=\{e,u,u^2,s,us,u^2s\}.$ The right adjoint action of $G=S_3.e$
on $M=e.S_3$ is given by
\[\ra \und{u}=(s,us,u^2s),\quad \ra \und{s} =(u,u^2)
(us,u^2s)\] while the left action of $e.S_3$ on $S_3.e$ is
trivial. The set $Z$ is $e.S_3$. It  splits into three conjugacy
classes $C^X$, namely $\{e\}$, $\{u,u^2\}$ and $\{s,us,u^2s\}$.

Following the general theory in Sections~\ref{classisection}, we
obtain in fact eight non-isomorphic irreducible bicovariant
differential calculi on $k(S_3)\lcross kS_3$ as follows:
$(i)$--$(ii)$ for $C^X=\{e\}$ we have one calculus of dimension 1
and one of dimension 2. $(iii)$--$(v)$ for $C^X=\{u,u^2\}$ we have
two calculi of dimension 2 and one of dimension 4.
$(iv)$--$(viii)$ for $C^X=\{s,us,u^2s\}$ we have two calculi of
dimension 3 and one of dimension 6. We omit details for most of
these calculi since they are similar in complexity and flavour to
Section~6.3, limiting ourselves to the most interesting one
$(viii)$ only:

\vspace{.5cm} $(viii)$\quad  $C^X=\{s,us,u^2s\},\;\;
\CM=<f_0,f_1,f_2,f_3,f_4,f_5>_k$,\quad where $q=e^{\frac{2\pi
i}{3}}$
\[f_0= e.s+q^2\und{u}.u^2s+q\und{u^2}.us ,\quad f_1=
\und{u}.s+q^2\und{u^2}.u^2s+qe.us\]
\[f_2=\und{u^2}.s+q^2e.u^2s+q\und{u}.us,\quad f_3=
\und{s}.s+q\und{u^2s}.us+q^2\und{us}.u^2s\]
\[f_4=\und{u^2s}.s+q\und{us}.us+q^2\und{s}.u^2s,\quad f_5=
\und{us}.s+q\und{s}.us+q^2\und{u^2s}.u^2s\] We have commutation
relations
\[e_if=R_{u^is}(f)e_i,\quad i\in \Z_6,\quad f\in k(S_3)\]
\[e_0s=se_3,\; e_3s=se_0,\quad e_is=se_{-i},\; \hbox{for}\; i\neq 0,3\]
\[e_2u=ue_0,\; e_5u=ue_3,\quad e_iu=ue_{i+1},\; \hbox{for}\; i\neq 2,5\]
and exterior differentials
\[\extd(f)=\partial_s(f)e_0+q\partial_{us}(f)e_1+q^2\partial_{u^2s}(f)e_2\]
\[ \theta=e_0+qe_1+q^2e_2,\quad \extd(u^i)=(q^{2i}-1)u^i\theta\]
\[\extd(u^is)=(q^{2i}-1)u^is\theta+q^{2i}s(e_3-e_0)+q^{2i+1}s(e_5-e_1)
+q^{2i+2}s(e_4-e_2),\quad i=0,1,2.\] The braiding is
\[\Psi(e_i\tens e_0)=qe_1\tens e_i,\quad \Psi(e_i\tens e_1)=q^2e_0\tens e_i,
\quad \Psi(e_i\tens e_2)=e_2\tens e_i,\]
\[\Psi(e_i\tens e_3)=q^2e_4\tens e_i,\quad \Psi(e_i\tens e_4)=qe_3\tens e_i,
\quad \Psi(e_i\tens e_5)=e_5\tens e_i,\quad i=2,5\]
\[\Psi(e_i\tens e_0)=q^2e_2\tens e_i,\quad \Psi(e_i\tens e_1)=e_1\tens e_i,
\quad \Psi(e_i\tens e_2)=qe_0\tens e_i,\]
\[\Psi(e_i\tens e_3)=qe_5\tens e_i,\quad \Psi(e_i\tens e_4)=e_4\tens e_i,
\quad \Psi(e_i\tens e_5)=q^2e_3\tens e_i,\quad i=1,4\]
\[\Psi(e_i\tens e_0)=e_0\tens e_i,\quad \Psi(e_i\tens e_1)=qe_2\tens e_i,
\quad \Psi(e_i\tens e_2)=q^2e_1\tens e_i,\]
\[\Psi(e_i\tens e_3)=e_3\tens e_i,\quad \Psi(e_i\tens e_4)=q^2e_5\tens e_i,
\quad \Psi(e_i\tens e_5)=qe_4\tens e_i,\quad i=0,3\] and yields
the degree 2 exterior algebra $\Omega^2(D(S_3))$ as 21-dimensional
with relations
\[e_i\wedge e_i=0,\quad i\in \Z_6,\quad \quad\{e_i,e_{i+3}\}=0,\quad i=0,1,2\]
\[e_0\wedge e_2+qe_2\wedge e_1+q^2e_1\wedge e_0=0,\quad
e_1\wedge e_2+qe_0\wedge e_1+q^2e_2\wedge e_0=0,\]
\[e_5\wedge e_3+qe_3\wedge e_4+q^2e_4\wedge e_5=0,\quad
e_5\wedge e_4+qe_3\wedge e_5+q^2e_4\wedge e_3=0\]
\[e_5\wedge e_0+e_3\wedge e_1+e_4\wedge e_2+q(e_0\wedge e_4+
e_1\wedge e_5+e_2\wedge e_3)=0\]
\[e_0\wedge e_5+e_1\wedge e_3+e_2\wedge e_4+q(e_4\wedge e_0+
e_5\wedge e_1+e_3\wedge e_2)=0.\] There are further relations in
degree 3, i.e. the entire Woronowicz exterior algebra in this
example is not quadratic. Its dimensions and cohomology in low
degree are
\[ \dim(\Omega)=1:6:21:60:152:\cdots,\quad H^0=k.1,\quad
H^1=k.\theta\oplus k.\bar\theta\] where
$\bar\theta=e_3+q^{-1}e_4+q^{-2}e_5$.\\

From these and similar computations for all the other calculi
(along the lines in Section~6.2) we find:
\begin{prop}
The 6-dimensional calculus $(viii)$ is the unique irreducible
calculus with $H^0=k.1$.
\end{prop}

We also have Poincar\'e duality at least where the exterior
algebra was small enough to be fully computed. For example, for
(vi)-(vii) the dimensions of $\Omega$ are 1:3:4:3:1 and the
dimensions of the cohomology are $6:6:0:6:6$.

The quantum double $D(S_3)$ is interesting for many reasons. Let
us note that being quasitriangular, it has a universal R-matrix or
quasitriangular structure $\mathcal{R}$ which controls the
noncocommutativity. This in turn is the nonAbelianness of the
underlying noncommutative group if one views $D(S_3)$ as a
function algebra, so should correspond to Riemannian curvature in
the setting of \cite{MAjRiema}. Our result is that there is a
unique irreducible calculus to take for this geometry. The ensuing
noncommutative Riemannian geometry will be developed elsewhere.

Also, from a mathematical point of view, $D(S_3)$ is a cotwist by
a multiplication-altering cocycle of the tensor product
$k(S_3)\tens kS_3$, its differential calculi can also  be obtained
from those of the tensor product $k(S_3)\tens kS_3$ by cotwisting
the exterior algebra according to the cotwisting theorem in
\cite{MAjOec}. This means that the classification of differential
calculi and their cohomology for $D(S_3)$ is exactly the same as
for the tensor product covered in Proposition~\ref{tensorcalc}.
Also note that until now the main example of a nontrivial
bicrossproduct in \cite{BM} was $k(\Z_6)\lrbicross k\Z_6$. We
find, however, that this is actually isomorphic as a Hopf algebra
to $(k(\Z_2)\rcrossco k\Z_3) \tens (k(\Z_3)\lcross k\Z_2)$ i.e. to
the tensor product $k(S_3)\tens kS_3$ again, hence has the same
features via twisting as $D(S_3)$. Similarly, replacing $\Z_6.S_3$
in Section~6.2 by the opposite factorization $S_3.\Z_6$ leads to
the dual bicrossproduct Hopf algebra $k(\Z_6)\lrbicross kS_3$
which by Proposition~2.1 in \cite{BM} is isomorphic to the quantum
double $D(S_3)$ again. Therefore several other known
bicrossroducts reduce to or have the same features as $D(S_3)$
above.

\subsubsection*{Acknowledgements} Main results and part of the writing
up was done when F.N was visiting Queen Mary, University of London
in the summer of 2002; he thanks the department there, D.~Lambert
from FUNDP Belgium and J-P. Antoine from UCL, Belgium
for organizing the visit for collaboration. \\

\end{document}